\RequirePackage{fix-cm}
\documentclass[smallextended]{svjour3}       
\smartqed  
\usepackage{graphicx}
\usepackage[colorlinks,linkcolor=red]{hyperref}   
\usepackage{bigstrut,multirow,rotating}      
\usepackage{booktabs}                             
\usepackage{enumerate}
\usepackage{amsmath}
\usepackage{bbm}
\usepackage{appendix}
\usepackage[numbers,sort&compress]{natbib}
\begin{document}

\title{Extracting Governing Laws from Sample Path Data of Non-Gaussian Stochastic Dynamical Systems
}

\titlerunning{Extracting Stochastic Governing Laws from Data}        

\author{Yang Li        \and
        Jinqiao Duan 
}


\institute{Li Y. \at
               School of Automation, Nanjing University of Science and Technology, Nanjing 210094, China           
           \and
           Duan J. \at
             Departments of Applied Mathematics \& Physics, Illinois Institute of Technology, Chicago, Illinois 60616, USA\\
\email{duan@iit.edu}
}

\date{Received: date / Accepted: date}

\maketitle

\begin{abstract}
Advances in  data science are leading to new progresses in the analysis and understanding of complex dynamics for systems with   experimental and observational data. With numerous physical phenomena exhibiting bursting, flights, hopping, and intermittent features, stochastic differential equations with non-Gaussian L\'evy noise are suitable to model these systems. Thus it is   desirable and essential to infer such equations from available data to reasonably predict dynamical behaviors. In this work, we consider a data-driven method to extract stochastic dynamical systems with non-Gaussian asymmetric (rather than the symmetric) L\'evy process, as well as Gaussian Brownian motion. We establish a theoretical framework and design a numerical algorithm to compute the asymmetric L\'evy jump measure,   drift and diffusion (i.e., nonlocal Kramers-Moyal formulas), hence obtaining the stochastic governing law,   from noisy data. Numerical experiments on several prototypical examples confirm the efficacy and accuracy of this method. This method will become an effective tool in discovering the governing laws from available data sets and in understanding the mechanisms underlying complex  random phenomena.
\end{abstract}

\section{Introduction}
\label{intro}

The extraction of physical laws from experimental or observable data is crucial to scientific and engineering applications where governing equations are unknown. Due to the deficiency of scientific understanding of mechanisms underlying complex phenomena, it is sometimes not feasible   to derive the explicit governing laws directly. These physical governing laws are often in the form of ordinary, partial or stochastic differential equations.

Therefore, there are recent   machine learning methods devoted to discovering governing laws of nonlinear phenomena from noisy data. Some of these works are based on Kramers-Moyal formula \cite{Ref8, DaiMinChaos, KM} or Koopman generator  \cite{Ref1, Ref2, WuFuDuan2019, LuYB2020, Yanxia2020}.  Other researchers have developed a data-driven method  called Sparse Identification of Nonlinear Dynamics, to learn ordinary \cite{Ref3}, partial \cite{Ref4, Ref5, Ref6, Ref7} and stochastic \cite{Ref8} differential equations from available data sets.

These techniques only focus on extracting either deterministic differential equations, or stochastic differential equations with Gaussian noise. However, there are   numerous systems involving random bursting, flights, intermittent, hopping, or rare transitions features in, for example, statistical physics \cite{Barn}, climate change \cite{Ref9}, gene regulation \cite{Ref10}, ecology \cite{Ref11}, and geophysical turbulence \cite{Ref12}. In consideration of their jump character, these systems are  suitable to be modeled as by stochastic differential equations with (non-Gaussian) L\'evy processes rather than with Gaussian fluctuations alone. For instance, Böttcher \cite{Ref13} presented a simple construction method for a class of stochastic processes   based on state space dependent mixing of  L\'evy processes. Ditlevsen revealed  that the climate change system may be modeled as stochastic differential equations with L\'evy process and Brownian motion, in  the context of the Greenland ice core measurement data \cite{Ref9}. Based on this assertion, Zheng et al. developed a probabilistic framework to investigate the maximum likelihood climate change for an energy balance system under the combined influence of greenhouse effect and non-Gaussian  $\alpha$-stable L\'evy motions \cite{Ref14}. Some researchers and we studied the escape phenomena \cite{Ref15, Ref16} and stochastic resonance \cite{Ref17, Ref18} of a neuron model driven by the non-Gaussian noise to detect its excitation behaviors. The non-Gaussian L\'evy motions are also used to characterize   random fluctuations in gene networks \cite{Ref19, Ref20, Ref21}, current-biased long Josephson junctions \cite{Ref22} and other scientific fields \cite{Ref23, Ref24}.

Recently, we devised a novel data-driven approach to extract stochastic governing equations with (idealized) symmetric  non-Gaussian L\'evy motion  from data \cite{Ref25}. In this present work, we  present this technique to the stochastic dynamical systems with  general asymmetric non-Gaussian L\'evy noise. Numerical experiments for prototypical examples verify its effectiveness.

This work is arranged as follows. In Section \ref{sec:2}, we derive formulas in the form of a theorem and a corollary to express the L\'evy jump measure, drift, and diffusion coefficient, in terms of either the transition probability density or sample paths. These may be regarded as nonlocal Kramers-Moyal formulas, in contrast to the usual (local) Kramers-Moyal formulas as in  \cite[Ch. 3]{KM}. Then we design a numerical algorithm to compute the jump measure, drift and diffusion coefficient, hence obtaining the stochastic governing law,  in Section \ref{sec:3}. We test our method by numerical experiments in Section \ref{sec:4}, and finally we conclude with Discussion in Section \ref{sec:5}.

\section{Method}
\label{sec:2}

In the previous section, we realize that random fluctuations often have both Gaussian and non-Gaussian statistical features.  By   L\'evy-It\^o decomposition theorem \cite{Ref29, Ref34},  a large class of  random fluctuations  are indeed modeled as linear combinations of  a (Gaussian) Brownian motion $\mathbf{B}_t$ and a (non-Gaussian) L\'evy process $\mathbf{L}_t$. We thus consider an $n$-dimensional stochastic dynamical system in the following form
\begin{equation} \label{eq1}
\textrm{d}\mathbf{x}\left( t \right)=\mathbf{b}\left( \mathbf{x}\left( t \right) \right)\textrm{d}t+\Lambda \left( \mathbf{x}\left( t \right) \right)\textrm{d}{{\mathbf{B}}_{t}}+\mathbf{\sigma }\textrm{d}{{\mathbf{L}}_{t}},
\end{equation}
where ${{\mathbf{B}}_{t}}={{\left[ {{B}_{1,t}},\ \cdots ,\ {{B}_{n,t}} \right]}^{T}}$  is an $n$-dimensional Brownian motion,  and ${{\mathbf{L}}_{t}}={{\left[ {{L}_{1,t}},\ \cdots ,\ {{L}_{n,t}} \right]}^{T}}$ is an $n$-dimensional non-Gaussian L\'evy process with independent components described in the Appendix \ref{App:A}. The vector $\mathbf{b}\left( \mathbf{x} \right)={{\left[ {{b}_{1}}\left( \mathbf{x} \right),\ \cdots ,\ {{b}_{n}}\left( \mathbf{x} \right) \right]}^{T}}$  is the drift coefficient (or vector field) in ${{\mathbbm{R}}^{n}}$  and $\Lambda \left( \mathbf{x} \right)$ is an $n\times n$  matrix.  The diffusion matrix is defined by $a\left( \mathbf{x} \right)=\Lambda {{\Lambda }^{T}}$.  We take $\mathbf{\sigma }$  as  a diagonal matrix with the positive diagonal component ${\sigma}_{i}$, indicating the noise intensity of the corresponding component of  L\'evy process. Assume that the initial condition is $\mathbf{x}\left( 0 \right)=\mathbf{z}$  and the jump measure of ${{L}_{i,t}}$ is ${{\nu }_{i}}\left( \textrm{d}\xi \right)={{W}_{i}}\left( \xi \right)\textrm{d}\xi$, with kernel $W$,  for $\xi\in \mathbbm{R}\backslash \left\{ 0 \right\}$.

\begin{figure}
\centering
  \includegraphics[height = 5cm, width = 6.5cm]{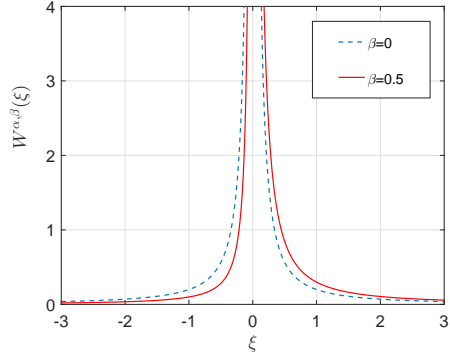}
\caption{Symmetric ($\beta=0$) and asymmetric ($\beta=0.5$) L\'evy kernel function ${{W}^{\alpha ,\beta }}\left( y \right)$ for $\alpha=0.5$.}
\label{fig:1}       
\end{figure}

For the sake of concreteness,  we consider a special but significant L\'evy process due to its extensive physical applications  \cite{Barn},  the so-called  $\alpha$-stable L\'evy motion \cite{Ref26, Ref27},     as an example to illustrate our method. Its detailed information is present at the Appendix \ref{App:A}. In this case, the  kernel function ${{W}^{\alpha ,\beta }}\left( \xi \right)$ has the following form
$$
{{W}^{\alpha ,\beta }}\left( \xi \right)=\left\{ \begin{array}{ccc}
   \frac{{{k}_{\alpha }}\left( 1+\beta  \right)}{2{{\left| \xi \right|}^{1+\alpha }}}, & \xi>0,  \\
   \frac{{{k}_{\alpha }}\left( 1-\beta  \right)}{2{{\left| \xi \right|}^{1+\alpha }}}, & \xi<0,  \\
\end{array} \right.
$$
where $\alpha$ is the stability parameter (or non-Gaussianity index) and
$$
{{k}_{\alpha }}=\left\{ \begin{array}{ccc}
   \frac{\alpha \left( 1-\alpha  \right)}{\Gamma \left( 2-\alpha  \right)\cos \left( {\pi \alpha }/{2}\; \right)}, & \alpha \ne 1,  \\
   \frac{2}{\pi }, & \alpha =1.  \\
\end{array} \right.
$$
The skewness parameter $\beta$ dominates the symmetry of  L\'evy motion, as shown in Fig. \ref{fig:1}.

According to Ref. \cite{Ref28, Ref29}, the Fokker-Planck equation for the probability density function $p\left( \mathbf{x},t|\mathbf{z},0 \right)$ for the solution of Eq. (\ref{eq1}) is
\begin{equation}\label{eq2}
\begin{aligned}
  \frac{\partial p}{\partial t}=& -\sum\limits_{i=1}^{n}{\frac{\partial }{\partial {{x}_{i}}}\left[ {{b}_{i}}p\left( \mathbf{x},t|\mathbf{z},0 \right) \right]}+\frac{1}{2}\sum\limits_{i,j=1}^{n}{\frac{{{\partial }^{2}}}{\partial {{x}_{i}}\partial {{x}_{j}}}\left[ {{a}_{ij}}p\left( \mathbf{x},t|\mathbf{z},0 \right) \right]} \\
 & -\sum\limits_{i=1}^{n}{\int_{\mathbbm{R}\backslash \left\{ 0 \right\}}{\left[ p\left( \mathbf{x},t|\mathbf{z},0 \right)-p\left( \mathbf{x}-{{\sigma }_{i}}{y}_{i}{{\mathbf{e}}_{i}},t|\mathbf{z},0 \right)-{{\sigma }_{i}}\chi _{i}^{\alpha }\left( {y}_{i} \right){y}_{i}\frac{\partial }{\partial {{x}_{i}}}p\left( \mathbf{x},t|\mathbf{z},0 \right) \right]W_{i}^{\alpha ,\beta }\left( {y}_{i} \right)\textrm{d}{y}_{i}}},
\end{aligned}
\end{equation}
where
$$
\chi _{i}^{\alpha }\left( {y}_{i} \right)=\left\{ \begin{matrix}
   0, & 0<\alpha <1,  \\
   {{\chi }_{\left| {y}_{i} \right|<\sigma _{i}^{-1}}}, & \alpha =1,  \\
   1, & 1<\alpha <2.  \\
\end{matrix} \right.
$$
Recall here that ${{\mathbf{e}}_{i}}$'s form the standard basis  in $\mathbbm{R}^n$,  and ${{\chi }_{\left| {y}_{i} \right|<\sigma _{i}^{-1}}}$ denotes the indicator function.  The initial condition is $p\left( \mathbf{x},0|\mathbf{z},0 \right)=\delta \left( \mathbf{x}-\mathbf{z} \right)$. Here the integral in the right hand side is understood as a Cauchy principal value integral.

In order to discover stochastic dynamical systems with non-Gaussian L\'evy noise as well as Gaussian noise from data, we derive the following theorem. It expresses jump measure, drift and diffusion in terms of  the solution of Fokker-Planck equation.

In what follows, the marginal probability distribution in $i$-th direction is denoted as ${{p}_{i}}\left( {{x}_{i}},t|{{z}_{i}},0 \right)$ and the set $\Gamma $ represents an $n$-dimensional cube ${{\left[ -\varepsilon ,\ \varepsilon  \right]}^{n}}$.

\newtheorem{thm}{\bf Theorem}
\begin{thm} (Relation between stochastic governing law and Fokker-Planck equation)\\
\label{thm1}
For every $\varepsilon >0$, the probability density function $p\left( \mathbf{x},t|\mathbf{z},0 \right)$ and the jump measure, drift and diffusion have the following relations:
\begin{enumerate}[1)]
\item For every ${{x}_{i}}$ and ${{z}_{i}}$ satisfying $\left| {{x}_{i}}-{{z}_{i}} \right|>\varepsilon $ and $i=1,\ 2,\ \ldots ,\ n$,
$$
\underset{t\to 0}{\mathop{\lim }}\,{{t}^{-1}}{{p}_{i}}\left( {{x}_{i}},t|{{z}_{i}},0 \right)=\sigma _{i}^{-1}W_{i}^{\alpha ,\beta }\left( \sigma _{i}^{-1}\left( {{x}_{i}}-{{z}_{i}} \right) \right)
$$
uniformly in ${{x}_{i}}$ and ${{z}_{i}}$.
\item For $i=1,\ 2,\ \ldots ,\ n$,
$$
\underset{t\to 0}{\mathop{\lim }}\,{{t}^{-1}}\int_{\mathbf{x}-\mathbf{z}\in \Gamma }{\left( {{x}_{i}}-{{z}_{i}} \right)p\left( \mathbf{x},t|\mathbf{z},0 \right) \textrm{d}\mathbf{x}}={{b}_{i}}\left( \mathbf{z} \right)+R_{i}^{\alpha ,\beta }\left( \varepsilon  \right),
$$
where $R_{i}^{\alpha ,\beta }\left( \varepsilon  \right)=\left\{ \begin{array}{ccc}
   \sigma _{i}^{-1}\int_{-\varepsilon }^{\varepsilon }{{y}_{i}W_{i}^{\alpha ,\beta }\left( \sigma _{i}^{-1}{y}_{i} \right)\textrm{d}{y}_{i}}, & \alpha <1,  \\
   \sigma _{i}^{-1}\left[ \int_{-\varepsilon }^{-1}{{y}_{i}W_{i}^{\alpha ,\beta }\left( \sigma _{i}^{-1}{y}_{i}\right)\textrm{d}{y}_{i}}+\int_{1}^{\varepsilon }{{y}_{i}W_{i}^{\alpha ,\beta }\left( \sigma _{i}^{-1}{y}_{i} \right)\textrm{d}{y}_{i}} \right], & \alpha =1,  \\
   -\sigma _{i}^{-1}\left[ \int_{-\infty }^{-\varepsilon }{{y}_{i}W_{i}^{\alpha ,\beta }\left( \sigma _{i}^{-1}{y}_{i} \right)\textrm{d}{y}_{i}}+\int_{\varepsilon }^{\infty }{{y}_{i}W_{i}^{\alpha ,\beta }\left( \sigma _{i}^{-1}{y}_{i} \right)\textrm{d}{y}_{i}} \right], & \alpha >1.  \\
\end{array} \right.$
\item For $i,j=1,\ 2,\ \ldots ,\ n$,
$$
\underset{t\to 0}{\mathop{\lim }}\,{{t}^{-1}}\int_{\mathbf{x}-\mathbf{z}\in \Gamma }{\left( {{x}_{i}}-{{z}_{i}} \right)\left( {{x}_{j}}-{{z}_{j}} \right)p\left( \mathbf{x},t|\mathbf{z},0 \right)\textrm{d}\mathbf{x}}={{a}_{ij}}\left( \mathbf{z} \right)+S_{ij}^{\alpha ,\beta }\left( \varepsilon  \right),
$$
where $S_{ii}^{\alpha ,\beta }\left( \varepsilon  \right)=\sigma _{i}^{-1}\int_{-\varepsilon }^{\varepsilon }{{{y}_{i}^{2}}W_{i}^{\alpha ,\beta }\left( \sigma _{i}^{-1}{y}_{i} \right)\textrm{d}{y}_{i}}$ and $S_{ij}^{\alpha ,\beta }\left( \varepsilon  \right)=0$ for $i\ne j$.
\end{enumerate}
\end{thm}

Note that for the purpose of the design of subsequent numerical scheme, we reformulate these relations in Theorem \ref{thm1} in the following corollary. Thus, the jump measure, drift and diffusion are expressed in terms of the sample paths of the stochastic differential equation (\ref{eq1}).

\newtheorem{cor}[thm]{\bf Corollary} 
\begin{cor}(Nonlocal Kramers-Moyal formulas)\\
\label{cor2}
For every $\varepsilon >0$, the sample path solution $\mathbf{x}\left( t \right)$ of the stochastic differential equation (\ref{eq1}) and the jump measure, drift and diffusion have the following relations:
\begin{enumerate}[1)]
\item For every ${{c}_{1}}$ and ${{c}_{2}}$ satisfying ${{c}_{1}}<{{c}_{2}}<0$ or $0<{{c}_{1}}<{{c}_{2}}$, and $i=1,\ 2,\ \ldots ,\ n$,
$$
\underset{t\to 0}{\mathop{\lim }}\,{{t}^{-1}}\mathbbm{P}\left\{ \left. {{x}_{i}}\left( t \right)-{{z}_{i}}\in \left[ {{c}_{1}},\ {{c}_{2}} \right) \right|\mathbf{x}\left( 0 \right)=\mathbf{z} \right\}=\sigma _{i}^{-1}\int_{{{c}_{1}}}^{{{c}_{2}}}{W_{i}^{\alpha ,\beta }\left( \sigma _{i}^{-1}{y}_{i} \right)\textrm{d}{y}_{i}};
$$
\item For $i=1,\ 2,\ \ldots ,\ n$,
$$
\begin{aligned}
  & \underset{t\to 0}{\mathop{\lim }}\,{{t}^{-1}}\mathbbm{P}\left\{ \left. \mathbf{x}\left( t \right)-\mathbf{z}\in \Gamma  \right|\mathbf{x}\left( 0 \right)=\mathbf{z} \right\}\cdot \text{E}\left[ \left. \left( {{x}_{i}}\left( t \right)-{{z}_{i}} \right) \right|\mathbf{x}\left( 0 \right)=\mathbf{z};\ \mathbf{x}\left( t \right)-\mathbf{z}\in \Gamma  \right] \\
 & ={{b}_{i}}\left( \mathbf{z} \right)+R_{i}^{\alpha ,\beta }\left( \varepsilon  \right); \\
\end{aligned}
$$
\item For $i,j=1,\ 2,\ \ldots ,\ n$,
$$
\begin{aligned}
  & \underset{t\to 0}{\mathop{\lim }}\,{{t}^{-1}}\mathbbm{P}\left\{ \left. \mathbf{x}\left( t \right)-\mathbf{z}\in \Gamma  \right|\mathbf{x}\left( 0 \right)=\mathbf{z} \right\}\cdot \text{E}\left[ \left. \left( {{x}_{i}}\left( t \right)-{{z}_{i}} \right)\left( {{x}_{j}}\left( t \right)-{{z}_{j}} \right) \right|\mathbf{x}\left( 0 \right)=\mathbf{z};\ \mathbf{x}\left( t \right)-\mathbf{z}\in \Gamma  \right] \\
 & ={{a}_{ij}}\left( \mathbf{z} \right)+S_{ij}^{\alpha ,\beta }\left( \varepsilon  \right). \\
\end{aligned}
$$
\end{enumerate}
\end{cor}

The three formulas in this corollary may be called  the nonlocal Kramers-Moyal formulas, in contrast to the usual (local) Kramers-Moyal formulas as in  \cite[Ch. 3]{KM}.  They express the jump measure, drift and diffusion   in terms of the sample paths of the stochastic differential equation (\ref{eq1}). In other words, with experimental or observational sample path data, we can extract the underlying stochastic differential equation.

The proofs of Theorem \ref{thm1} and Corollary \ref{cor2} are   in Appendix \ref{App:B} and \ref{App:C}, respectively.

\section{Numerical algorithm}
\label{sec:3}

We now devise a numerical algorithm to  extract a   stochastic differential equation model from its sample path data.
Assume that there exists a pair of data sets containing $M$ elements, respectively,
\begin{equation}\label{eq3}
\begin{aligned}
  & Z=\left[ {{\mathbf{z}}_{1}},\ {{\mathbf{z}}_{2}},\ \cdots ,\ {{\mathbf{z}}_{M}} \right], \\
 & X=\left[ {{\mathbf{x}}_{1}},\ {{\mathbf{x}}_{2}},\ \cdots ,\ {{\mathbf{x}}_{M}} \right], \\
\end{aligned}
\end{equation}
where ${{\mathbf{x}}_{j}}$ is the image of ${{\mathbf{z}}_{j}}$ after a small evolution time $h$. In addition, we also choose a dictionary of basis functions $\Psi \left( \mathbf{x} \right)=\left[ {{\psi }_{1}}\left( \mathbf{x} \right),\ {{\psi }_{2}}\left( \mathbf{x} \right),\ \cdots ,\ {{\psi }_{K}}\left( \mathbf{x} \right) \right]$. With the data sets, the dictionary and the Corollary \ref{cor2} in the previous section, we now propose the following algorithm to identify the kernel function and noise intensity of the L\'evy motion, the drift coefficient and the diffusion matrix.

\subsection{Algorithm for identification of the L\'evy motion}
\label{subsec:3.1}
In order to determine the L\'evy motion, we need to identify all the parameters ${{\alpha }_{i}}$, ${{\beta }_{i}}$ and ${{\sigma }_{i}}$ for $i=1,\ 2,\ \ldots ,\ n$.  From the first assertion of Corollary \ref{cor2},  we see that the two sides of the equation only depend on the difference of ${{x}_{i}}$ and ${{z}_{i}}$ rather than their specific positions. Thus we construct a new data set ${{Y}_{i}}=\left[ y_{i}^{1},\ y_{i}^{2},\ \cdots ,\ y_{i}^{M} \right]$ with ${{y}_{i}}={{x}_{i}}-{{z}_{i}}$. Therefore, the probability in the left hand side of the first assertion can be approximated by the ratio of the number of the points falling into the interval $\left[ {{c}_{1}},\ {{c}_{2}} \right)$ to the total number $M$ in ${{Y}_{i}}$.

\begin{figure}
\centering
  \includegraphics[height = 5cm, width = 6.5cm]{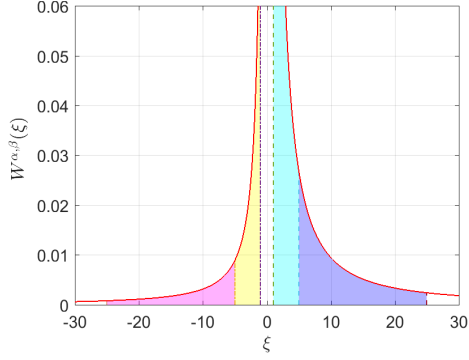}
\caption{Integration intervals for $\varepsilon =1$, $m=5$ and $N=1$.}
\label{fig:2}       
\end{figure}

Specifically, we consider the integral on $2N+2$ intervals $\left[ \varepsilon ,\ m\varepsilon  \right)$, $\left[ m\varepsilon ,\ {{m}^{2}}\varepsilon  \right)$, $\ldots$, $\left[ {{m}^{N}}\varepsilon ,\ {{m}^{N+1}}\varepsilon  \right)$ and $\left[ -m\varepsilon ,\ -\varepsilon  \right)$, $\left[ -{{m}^{2}}\varepsilon ,\ -m\varepsilon  \right)$, $\ldots$, $\left[ -{{m}^{N+1}}\varepsilon ,\ -{{m}^{N}}\varepsilon  \right)$ with the positive integer $N$, the positive real number $\varepsilon $ and the real number $m>1$. This is  illustrated in Fig. \ref{fig:2}. For the convenience of representations, we neglect the subscript “$i$” of ${{\alpha }_{i}}$, ${{\beta }_{i}}$ and ${{\sigma }_{i}}$ in the following equations. Assume that there are $n_{0}^{+},\ n_{1}^{+},\ \ldots ,\ n_{N}^{+}$ and $n_{0}^{-},\ n_{1}^{-},\ \ldots ,\ n_{N}^{-}$ points from the data set ${{Y}_{i}}$ which fall into these intervals respectively. Therefore,
$$
{{h}^{-1}}\mathbbm{P}\left\{ \left. {{x}_{i}}\left( h \right)-{{z}_{i}}\in \left[ {{m}^{k}}\varepsilon ,\ {{m}^{k+1}}\varepsilon  \right) \right|\mathbf{x}\left( 0 \right)=\mathbf{z} \right\}\approx {{h}^{-1}}{{M}^{-1}}n_{k}^{+}.
$$
On the other hand, the integration from the right hand side yields
$$
\begin{aligned}
  & {{\sigma }^{-1}}\int_{{{m}^{k}}\varepsilon }^{{{m}^{k+1}}\varepsilon }{W_{i}^{\alpha ,\beta }\left( {{\sigma }^{-1}}{y}_{i} \right)\textrm{d}{y}_{i}} \\
 & ={{{\sigma }^{\alpha }}{{k}_{\alpha }}}\frac{\left( 1+\beta  \right)}{2}\int_{{{m}^{k}}\varepsilon }^{{{m}^{k+1}}\varepsilon }{{{\left| {y}_{i} \right|}^{-\left( 1+\alpha  \right)}}\text{d}{y}_{i}} \\
 & ={{{\sigma }^{\alpha }}{{k}_{\alpha }}\left( 1+\beta  \right){{\alpha }^{-1}}{{\varepsilon }^{-\alpha }}{{m}^{-k\alpha }}\left( 1-{{m}^{-\alpha }} \right)}/{2}\;. \\
\end{aligned}
$$
Combining the two equations, we have $N+1$ equalities
\begin{equation}\label{eq4}
{{{\sigma }^{\alpha }}{{k}_{\alpha }}\left( 1+\beta  \right){{\alpha }^{-1}}{{\varepsilon }^{-\alpha }}{{m}^{-k\alpha }}\left( 1-{{m}^{-\alpha }} \right)}/{2}={{h}^{-1}}{{M}^{-1}}n_{k}^{+},\quad k=0,\ 1,\ \ldots ,\ N.
\end{equation}
On the other $N+1$ intervals, we also have
\begin{equation}\label{eq5}
{{{\sigma }^{\alpha }}{{k}_{\alpha }}\left( 1-\beta  \right){{\alpha }^{-1}}{{\varepsilon }^{-\alpha }}{{m}^{-k\alpha }}\left( 1-{{m}^{-\alpha }} \right)}/{2}={{h}^{-1}}{{M}^{-1}}n_{k}^{-},\quad k=0,\ 1,\ \ldots ,\ N.
\end{equation}
The summations of Eqs. (\ref{eq4}) and (\ref{eq5}) yield
$$
{{{\sigma }^{\alpha }}{{k}_{\alpha }}{{\alpha }^{-1}}{{\varepsilon }^{-\alpha }}{{m}^{-k\alpha }}\left( 1-{{m}^{-\alpha }} \right)}={{h}^{-1}}{{M}^{-1}}(n_{k}^{+}+n_{k}^{-}),\quad k=0,\ 1,\ \ldots ,\ N.
$$
The ratios of the first equation ($k=0$) to the other $N$ equations lead to the solutions
\begin{equation}\label{eq6}
\alpha ={{\left( k\ln m \right)}^{-1}}\ln \frac{n_{0}^{+}+n_{0}^{-}}{n_{k}^{+}+n_{k}^{-}},\quad k=1,\ 2,\ \ldots ,\ N.
\end{equation}
If $N=1$, it is the optimal solution $\tilde{\alpha }$ of the parameter $\alpha$. In order to make full use of the data information, a bigger $N$ can be selected to identify the optimal solution $\tilde{\alpha }$ as the mean value of Eqs. (\ref{eq6}). Let $\rho ={\sum\limits_{k}{n_{k}^{-}}}/{\sum\limits_{k}{n_{k}^{+}}}$. Then the ratio of the summation of Eqs. (\ref{eq4}) over the summation of Eqs. (\ref{eq5}) yields
\begin{equation}\label{eq7}
\beta =\frac{1-\rho }{1+\rho }.
\end{equation}
Therefore, the noise intensity $\sigma$ is computed as
\begin{equation}\label{eq8}
\sigma ={{\left[ \frac{\tilde{\alpha }{{\varepsilon }^{{\tilde{\alpha }}}}{{m}^{k\tilde{\alpha }}}\left( n_{k}^{+}+n_{k}^{-} \right)}{{{k}_{{\tilde{\alpha }}}}hM\left( 1-{{m}^{-\tilde{\alpha }}} \right)} \right]}^{{1}/{{\tilde{\alpha }}}\;}},\quad k=0,\ 1,\ \ldots ,\ N.
\end{equation}
Hence, the optimal noise intensity $\tilde{\sigma }$ can be identified as the mean value of Eqs. (\ref{eq8}). By going through all $n$ dimensions for $i$, we obtain $3n$ parameters ${{\alpha }_{i}}$, ${{\beta }_{i}}$ and ${{\sigma }_{i}}$ to determine the L\'evy jump measure as well as its noise intensity.

Note that if we consider the symmetric  L\'evy motion, then $\rho =1$ and $\beta =0$ in Eq. (\ref{eq7}). Let ${{n}_{k}}=n_{k}^{+}+n_{k}^{-}$, thus Eqs. (\ref{eq6}) and (\ref{eq8}) recover the computations of the symmetric case. It is also important to emphasize that this method avoids the problem of curse of dimensionality since we can deal with every component independently. In other words, the increase of dimensionality does not necessarily require more data.

\subsection{Algorithm for identification of the drift term}
\label{subsec:3.2}
In order to identify the drift coefficient $\mathbf{b}\left( \mathbf{x} \right)$, we approximate its every component in terms of the dictionary of basis functions $\Psi \left( \mathbf{x} \right)$ as ${{b}_{i}}\left( \mathbf{x} \right)\approx \sum\limits_{k=1}^{K}{{{c}_{i,k}}{{\psi }_{k}}\left( \mathbf{x} \right)}$, $i=1,\ 2,\ \ldots ,\ n$. According to the second assertion of Corollary \ref{cor2}, the computation of the drift term only requires the data within the cube $\mathbf{x}-\mathbf{z}\in \Gamma $. Therefore, after deleting the data outside the cube $\Gamma$ in the data sets $Z$ and $M$, we obtain the new data sets with $\hat{M}$ elements respectively
\begin{equation}\label{eq9}
\begin{aligned}
  & \hat{Z}=\left[ {{{\mathbf{\hat{z}}}}_{1}},\ {{{\mathbf{\hat{z}}}}_{2}},\ \cdots ,\ {{{\mathbf{\hat{z}}}}_{{\hat{M}}}} \right], \\
 & \hat{X}=\left[ {{{\mathbf{\hat{x}}}}_{1}},\ {{{\mathbf{\hat{x}}}}_{2}},\ \cdots ,\ {{{\mathbf{\hat{x}}}}_{{\hat{M}}}} \right]. \\
\end{aligned}
\end{equation}
Based on the second assertion of Corollary \ref{cor2}, we derive the following group of equations
\begin{equation}\label{eq10}
\begin{aligned}
  & A{{\mathbf{c}}_{i}}={{B}_{i}}, \\
 & A=\left[ \begin{array}{ccc}
   {{\psi }_{1}}\left( {{{\mathbf{\hat{z}}}}_{1}} \right) & \cdots  & {{\psi }_{K}}\left( {{{\mathbf{\hat{z}}}}_{1}} \right)  \\
   \vdots  & \ddots  & \vdots   \\
   {{\psi }_{1}}\left( {{{\mathbf{\hat{z}}}}_{{\hat{M}}}} \right) & \cdots  & {{\psi }_{K}}\left( {{{\mathbf{\hat{z}}}}_{{\hat{M}}}} \right)  \\
\end{array} \right], \\
 & {{B}_{i}}=\hat{M}{{M}^{-1}}{{h}^{-1}}{{\left[ {{{\hat{x}}}_{i,1}}-{{{\hat{z}}}_{i,1}},\ \cdots ,\ {{{\hat{x}}}_{i,\hat{M}}}-{{{\hat{z}}}_{i,\hat{M}}} \right]}^{T}}-R_{i}^{\alpha ,\beta }\left( \varepsilon  \right). \\
\end{aligned}
\end{equation}
Here, the computation of $R_{i}^{\alpha ,\beta }\left( \varepsilon  \right)$ requires the estimated parameters ${{\alpha }_{i}}$, ${{\beta }_{i}}$ and ${{\sigma }_{i}}$ in Subsection \ref{subsec:3.1}. This group of equations may be solved via the least squares approximation $\underset{{{\mathbf{c}}_{i}}}{\mathop{\min }}\,{{\left\| A{{\mathbf{c}}_{i}}-{{B}_{i}} \right\|}_{2}^{2}}$. It is well known that its optimal solution is provided by
\begin{equation}\label{eq11}
{{\mathbf{\tilde{c}}}_{i}}={{\left( {{A}^{T}}A \right)}^{-1}}\left( {{A}^{T}}{{B}_{i}} \right).
\end{equation}

\subsection{Algorithm for identification of the diffusion term}
\label{subsec:3.3}
Finally, the diffusion matrix $a\left( \mathbf{x} \right)$ remains to be found. Based on the dictionary of basis functions $\Psi \left( \mathbf{x} \right)$, its component can be estimated as ${{a}_{ij}}\left( \mathbf{x} \right)\approx \sum\limits_{k=1}^{K}{{{d}_{ij,k}}{{\psi }_{k}}\left( \mathbf{x} \right)}$, $i,j=1,\ 2,\ \ldots ,\ n$. According to the third assertion of Corollary \ref{cor2}, we derive the following group of equations
\begin{equation}\label{eq12}
\begin{aligned}
  & A{{\mathbf{d}}_{ij}}={{B}_{ij}}, \\
 & A=\left[ \begin{array}{ccc}
   {{\psi }_{1}}\left( {{{\mathbf{\hat{z}}}}_{1}} \right) & \cdots  & {{\psi }_{K}}\left( {{{\mathbf{\hat{z}}}}_{1}} \right)  \\
   \vdots  & \ddots  & \vdots   \\
   {{\psi }_{1}}\left( {{{\mathbf{\hat{z}}}}_{{\hat{M}}}} \right) & \cdots  & {{\psi }_{K}}\left( {{{\mathbf{\hat{z}}}}_{{\hat{M}}}} \right)  \\
\end{array} \right], \\
 & {{B}_{ij}}=\hat{M}{{M}^{-1}}{{h}^{-1}}{{\left[ \left( {{{\hat{x}}}_{i,1}}-{{{\hat{z}}}_{i,1}} \right)\left( {{{\hat{x}}}_{j,1}}-{{{\hat{z}}}_{j,1}} \right),\ \cdots ,\ \left( {{{\hat{x}}}_{i,\hat{M}}}-{{{\hat{z}}}_{i,\hat{M}}} \right)\left( {{{\hat{x}}}_{j,\hat{M}}}-{{{\hat{z}}}_{j,\hat{M}}} \right) \right]}^{T}}-S_{ij}^{\alpha ,\beta }\left( \varepsilon  \right). \\
\end{aligned}
\end{equation}
Subsequently, we also consider the least squares problem $\underset{{{\mathbf{d}}_{ij}}}{\mathop{\min }}\,{{\left\| A{{\mathbf{d}}_{ij}}-{{B}_{ij}} \right\|}_{2}^{2}}$ to find the optimal solution. It is well known that this leads to the solution
\begin{equation}\label{eq13}
{{\mathbf{\tilde{d}}}_{ij}}={{\left( {{A}^{T}}A \right)}^{-1}}\left( {{A}^{T}}{{B}_{ij}} \right).
\end{equation}
Since the diffusion matrix is symmetric, we just need to compute the coefficients for $i=1,\ 2,\ \ldots ,\ n$, $j=i,\ i+1,\ \ldots ,\ n$.

Remark that it is not unique to infer the coefficient $\Lambda \left( \mathbf{x} \right)$ in terms of the diffusion matrix $a\left( \mathbf{x} \right)$. Via the symmetry of $a\left( \mathbf{x} \right)$, the orthogonal matrix $Q$ that is composed of the eigenvectors of $a\left( \mathbf{x} \right)$ can diagonalize the diffusion matrix, i.e., $J={{Q}^{T}}aQ$ with the component in the diagonal ${{J}_{ii}}\ge 0$. Then $J=\left( {{Q}^{T}}\Lambda  \right)\cdot {{\left( {{Q}^{T}}\Lambda  \right)}^{T}}$. Denoting $\sqrt{J}$ as the diagonal matrix with the diagonal component $\sqrt{{{J}_{ii}}}$,   the matrix $\Lambda =Q\sqrt{J}$ is a solution of this equation. Note that for arbitrary orthogonal matrix $\Lambda =Q\sqrt{J}$, the matrix $Q\sqrt{J}T$ is also a solution. Even so, they are statistically equivalent as they have the same Fokker-Planck equation for a given diffusion matrix.

\section{Examples}
\label{sec:4}
We now present two prototypical examples to demonstrate our method for discovering stochastic dynamical systems from sample path data sets.

\newtheorem{exa}{\bf Example}
\begin{exa}\label{exa1}
Consider a  stochastic Lorenz   system
$$
\begin{aligned}
  & \textrm{d}{{x}_{1}}=\textrm{10}\left( -{{x}_{1}}+{{x}_{2}} \right)\textrm{d}t+\left( 1+{{x}_{3}} \right)\textrm{d}{{B}_{1,t}}+\textrm{d}{{B}_{2,t}}+2\textrm{d}{{L}_{1,t}}, \\
 & \textrm{d}{{x}_{2}}=\left( 4{{x}_{1}}-{{x}_{2}}-{{x}_{1}}{{x}_{3}} \right)\textrm{d}t+{{x}_{2}}\textrm{d}{{B}_{2,t}}+\textrm{d}{{L}_{2,t}}, \\
 & \textrm{d}{{x}_{3}}=\left( -{8}/{3}\;{{x}_{3}}+{{x}_{1}}{{x}_{2}} \right)\textrm{d}t+{{x}_{1}}\textrm{d}{{B}_{3,t}}+0.5\textrm{d}{{L}_{3,t}}. \\
\end{aligned}
$$
The L\'evy noise intensity ${{\sigma }_{1}}=2$, ${{\sigma }_{2}}=1$, ${{\sigma }_{3}}=0.5$ and the drift and the diffusion coefficients
$$
\begin{aligned}
  & \mathbf{b}\left( \mathbf{x} \right)={{\left[ 10\left( -{{x}_{1}}+{{x}_{2}} \right),\ 4{{x}_{1}}-{{x}_{2}}-{{x}_{1}}{{x}_{3}},\ -\frac{8}{3}{{x}_{3}}+{{x}_{1}}{{x}_{2}} \right]}^{T}}, \\
 & a\left( \mathbf{x} \right)=\left[ \begin{array}{cccc}
   2+2{{x}_{3}}+x_{3}^{2} & \ &{{x}_{2}} & 0  \\
   {{x}_{2}} & \ & x_{2}^{2} & 0  \\
   0 & \ &  0 & x_{1}^{2}  \\
\end{array} \right]. \\
\end{aligned}
$$
The time step is fixed as $h=0.001$ and the chosen initial points $Z=\left[ {{\mathbf{z}}_{1}}, {{\mathbf{z}}_{2}}, \ldots , {{\mathbf{z}}_{M}} \right]$ are distributed uniformly within the  cube  $\left[ -2,2 \right]\times \left[ -2, 2 \right]\times \left[ -2, 2 \right]$ with a mesh $400\times 400\times 400$.   The data set $X$ is generated by Euler scheme of numerical scheme starting  from $Z$ after the time $h$ with ${{\alpha }_{1}}=0.5$, ${{\beta }_{1}}=0.5$, ${{\alpha }_{2}}=1$, ${{\beta }_{2}}=0$ and ${{\alpha }_{3}}=1.5$, ${{\beta }_{3}}=-0.5$. The dictionary is selected as polynomials to estimate the drift and diffusion terms.

With the parameters fixed as $N=2$, $\varepsilon =1$ and $m=5$, we compute the three groups of the parameters $\alpha$, $\beta$ and $\sigma$, and the drift and diffusion terms via our  proposed algorithms and list the results in Tables \ref{tab:1}, \ref{tab:2} and \ref{tab:3}. As we   see, the results agree well with the true parameters.

\begin{table}[htbp]
  \centering
  \caption{IdentifiedL\'evy motion for the three-dimensional Lorenz system}
    \begin{tabular}{ccccccc}
    \toprule
    \multicolumn{1}{c}{\multirow{2}[2]{*}{Parameter}} & \multicolumn{2}{c}{${L}_{1}$} & \multicolumn{2}{c}{${L}_{2}$} & \multicolumn{2}{c}{${L}_{3}$} \\
          & True  & \multicolumn{1}{p{4.04em}}{Learned} & True  & \multicolumn{1}{p{4.04em}}{Learned} & True  & \multicolumn{1}{p{4.04em}}{Learned} \\
    \midrule
    $\alpha$      & 0.5   & 0.498 & 1     & 1.0038 & 1.5   & 1.4832 \\
    $\beta$      & 0.5   & 0.5066 & 0     & 0.0003 & -0.5  & -0.5118 \\
    $\sigma$      & 2     & 2.01  & 1     & 0.9956 & 0.5   & 0.4907 \\
    \bottomrule
    \end{tabular}%
  \label{tab:1}%
\end{table}%

\begin{table}[htbp]
  \centering
  \caption{Identified drift term for the three-dimensional Lorenz system}
    \begin{tabular}{ccccccc}
    \toprule
    \multicolumn{1}{c}{\multirow{2}[1]{*}{Basis}} & \multicolumn{2}{c}{${{b}_{1}}\left( \mathbf{x} \right)$} & \multicolumn{2}{c}{${{b}_{2}}\left( \mathbf{x} \right)$} & \multicolumn{2}{c}{${{b}_{3}}\left( \mathbf{x} \right)$} \\
          & True  & \multicolumn{1}{p{4.04em}}{Learned} & True  & \multicolumn{1}{p{4.04em}}{Learned} & True  & \multicolumn{1}{p{4.04em}}{Learned} \\
    \midrule
    \multicolumn{1}{c}{1} & 0     & 0     & 0     & 0     & 0     & 0 \\
    \multicolumn{1}{c}{${x}_{1}$}      & -10   & -9.9801 & 4     & 3.9941 & 0     & 0 \\
    \multicolumn{1}{c}{${x}_{2}$}      & 10    & 9.9666 & -1    & -1.0036 & 0     & 0 \\
    \multicolumn{1}{c}{${x}_{3}$}      & 0     & 0     & 0     & 0     & -2.6667 & -2.6502 \\
    \multicolumn{1}{c}{$x_{1}^{2}$}      & 0     & 0     & 0     & 0     & 0     & 0 \\
    \multicolumn{1}{c}{${x}_{1}{x}_{2}$}      & 0     & 0     & 0     & 0     & 1     & 1 \\
    \multicolumn{1}{c}{${x}_{1}{x}_{3}$}      & 0     & 0     & -1    & -0.9917 & 0     & 0 \\
    \multicolumn{1}{c}{$x_{2}^{2}$}      & 0     & 0     & 0     & 0     & 0     & 0 \\
    \multicolumn{1}{c}{${x}_{2}{x}_{3}$}      & 0     & 0     & 0     & 0     & 0     & 0 \\
    \multicolumn{1}{c}{$x_{3}^{2}$}      & 0     & 0     & 0     & 0     & 0     & 0 \\
    \bottomrule
    \end{tabular}%
  \label{tab:2}%
\end{table}%

\begin{table}[htbp]
  \centering
  \caption{Identified diffusion term for the three-dimensional Lorenz system}
    \begin{tabular}{ccccccccccccc}
    \toprule
    \multicolumn{1}{c}{\multirow{2}[1]{*}{Basis}} & \multicolumn{2}{c}{${{a}_{11}}\left( \mathbf{x} \right)$} & \multicolumn{2}{c}{${{a}_{12}}\left( \mathbf{x} \right)$} & \multicolumn{2}{c}{${{a}_{13}}\left( \mathbf{x} \right)$} & \multicolumn{2}{c}{${{a}_{22}}\left( \mathbf{x} \right)$} & \multicolumn{2}{c}{${{a}_{23}}\left( \mathbf{x} \right)$} & \multicolumn{2}{c}{${{a}_{33}}\left( \mathbf{x} \right)$} \\
          & True  & \multicolumn{1}{p{4.04em}}{Learned} & True  & \multicolumn{1}{p{4.04em}}{Learned} & True  & \multicolumn{1}{p{4.04em}}{Learned} & True  & \multicolumn{1}{p{4.04em}}{Learned} & True  & \multicolumn{1}{p{4.04em}}{Learned} & True  & \multicolumn{1}{p{4.04em}}{Learned} \\
    \midrule
    \multicolumn{1}{c}{1} & 2     & 2.2748 & 0     & 0     & 0     & 0     & 0     & 0     & 0     & 0     & 0     & 0 \\
     \multicolumn{1}{c}{${x}_{1}$}      & 0     & 0     & 0     & 0     & 0     & 0     & 0     & 0     & 0     & 0     & 0     & 0 \\
     \multicolumn{1}{c}{${x}_{2}$}      & 0     & 0     & 1     & 0.9968 & 0     & 0     & 0     & 0     & 0     & 0     & 0     & 0 \\
    \multicolumn{1}{c}{${x}_{3}$}       & 2     & 1.993 & 0     & 0     & 0     & 0     & 0     & 0     & 0     & 0     & 0     & 0 \\
   \multicolumn{1}{c}{$x_{1}^{2}$}       & 0     & 0     & 0     & 0     & 0     & 0     & 0     & 0     & 0     & 0     & 1     & 1.0137 \\
  \multicolumn{1}{c}{${x}_{1}{x}_{2}$}        & 0     & 0     & 0     & 0     & 0     & 0     & 0     & 0     & 0     & 0     & 0     & 0 \\
  \multicolumn{1}{c}{${x}_{1}{x}_{3}$}        & 0     & 0     & 0     & 0     & 0     & 0     & 0     & 0     & 0     & 0     & 0     & 0 \\
  \multicolumn{1}{c}{$x_{2}^{2}$}        & 0     & 0     & 0     & 0     & 0     & 0     & 1     & 1.0139 & 0     & 0     & 0     & 0 \\
  \multicolumn{1}{c}{${x}_{2}{x}_{3}$}        & 0     & 0     & 0     & 0     & 0     & 0     & 0     & 0     & 0     & 0     & 0     & 0 \\
 \multicolumn{1}{c}{$x_{3}^{2}$}         & 1     & 0.9982 & 0     & 0     & 0     & 0     & 0     & 0     & 0     & 0     & 0     & 0 \\
    \bottomrule
    \end{tabular}%
  \label{tab:3}%
\end{table}%
\end{exa}

\begin{exa}\label{exa2}
Consider a gene regulation  system driven by both Gaussian Brownian noise and non-Gaussian L\'evy noise with a rational drift coefficient \cite{Ref30}
$$
\textrm{d}x\left( t \right)=\left[ \frac{{{k}_{f}}{{x}^{2}}\left( t \right)}{{{x}^{2}}\left( t \right)+{{K}_{d}}}-{{k}_{d}}x\left( t \right)+{{R}_{bas}} \right]\textrm{d}t+\frac{x\left( t \right)}{\sqrt{{{x}^{2}}\left( t \right)+0.5}}\textrm{d}{{B}_{t}}+0.5\textrm{d}{{L}_{t}},
$$
where the system parameters are ${{k}_{f}}=6 {{\min }^{-1}}$, ${{K}_{d}}=10$, ${{k}_{d}}=1 {{\min }^{-1}}$, and ${{R}_{bas}}=0.4 {{\min }^{-1}}$. Then we have the drift coefficient $b\left( x \right)={{{k}_{f}}{{x}^{2}}}/{\left( {{x}^{2}}+{{K}_{d}} \right)}-{{k}_{d}}x+{{R}_{bas}}$, the diffusion coefficient $a\left( x \right)={{{x}^{2}}}/{\left( {{x}^{2}}+0.5 \right)}$ and the L\'evy noise intensity $\sigma =0.5$. The time step is fixed as $h=0.001$ and the chosen ${{10}^{7}}$ initial points $Z=\left[ {{z}_{1}},\ {{z}_{2}},\ \ldots ,\ {{z}_{M}} \right]$ are distributed uniformly within the interval $\left[ 0,\ 5 \right]$. The simulated data set $X$ is generated by Euler scheme of numerical scheme  starting from $Z$ after the time $h$. In order to check out the validation of our method to non-polynomial basis functions, we choose the following dictionary as in Ref. \cite{Ref8}
$$
\begin{aligned}
  \Psi \left( x \right)=& \left[ 1,\ x,\ {{x}^{2}},\ {{x}^{3}},\ \sin x,\ \cos 11x,\ \sin 11x, \right. \\
 & -10{{\tanh }^{2}}\left( 10x \right)+10,\ -10{{\tanh }^{2}}\left( 10x-10 \right)+10, \\
 & \exp \left\{ -50{{x}^{2}} \right\},\ \exp \left\{ -50{{\left( x-3 \right)}^{2}} \right\},\ \exp \left\{ -0.3{{x}^{2}} \right\}, \\
 & \exp \left\{ -0.3{{\left( x-3 \right)}^{2}} \right\},\ \exp \left\{ -2{{\left( x-2 \right)}^{2}} \right\},\ \exp \left\{ -50{{\left( x-4 \right)}^{2}} \right\}, \\
 & \exp \left\{ -0.6{{\left( x-4 \right)}^{2}} \right\},\ \exp \left\{ -0.6{{\left( x-3 \right)}^{2}} \right\}, \\
 & \left. -2{{\tanh }^{2}}\left( 2x-4 \right)+2,\ {{\tanh }^{2}}\left( x-4 \right)+1 \right].
\end{aligned}
$$

Here we only show the  case with  $\alpha =1.5$ and $\beta =-0.5$.   We choose the parameters as $N=2$, $\varepsilon =1$ and $m=5$. According to our  proposed method in Section \ref{sec:3}, the learned values of the stability parameter $\alpha$, the skewness parameter $\beta$ and the L\'evy noise intensity $\sigma$ are listed in Table \ref{tab:4}. It is seen that they are consistent with the true values. Employing Eqs. (\ref{eq11}) and (\ref{eq13}), we compute the least square solutions $\mathbf{\tilde{c}}$ and $\mathbf{\tilde{d}}$ as
$$
\begin{aligned}
  \mathbf{\tilde{c}}=& \left[ -118.8168,\ -28.0260,\ 18.4281,\ -2.0028, \right. \\
 & 20.5302,\ 0,\ 0,\ 0,\ 0,\ 0,\ 0,\ 62.0744,\ 50.5969, \\
 & {{\left. 0,\ 6.2844,\ 0,\ 20.0348,\ -20.7253,\ 0,\ 26.9544 \right]}^{T}}, \\
 \mathbf{\tilde{d}}=& \left[ -121.6943,\ -30.3672,\ 20.5902,\ -2.2510, \right. \\
 & 24.4084,\ 0,\ 0,\ 0,\ 0,\ 0,\ 0,\ 68.1945,49.7766, \\
 & {{\left. 0,5.9156,\ 0,19.6789,-19.9246,0,\ 25.1716 \right]}^{T}}.
\end{aligned}
$$
According to these learned coefficients, the estimated drift and diffusion terms are plotted in Fig. \ref{fig:3}, indicating that the learned results   provide a  reasonable  approximation in parameter range of interest. This implies that our method is valid even for stochastic dynamical systems with    non-polynomial drift.

\begin{table}[htbp]
  \centering
  \caption{IdentifiedL\'evy motion for the genetic regulatory system}
    \begin{tabular}{ccc}
    \toprule
    \multicolumn{1}{p{4.04em}}{Parameter} & True  & \multicolumn{1}{p{4.04em}}{Learned} \\
    \midrule
     \multicolumn{1}{c}{$\alpha$}     & 1.5   & 1.5406 \\
     \multicolumn{1}{c}{$\beta$}     & -0.5  & -0.5237 \\
     \multicolumn{1}{c}{$\sigma$}     & 0.5   & 0.538 \\
    \bottomrule
    \end{tabular}%
  \label{tab:4}%
\end{table}%

\begin{figure}
	\begin{minipage}{0.48\linewidth}
		\leftline{(a)}
		\centerline{\includegraphics[height = 5cm, width = 6.5cm]{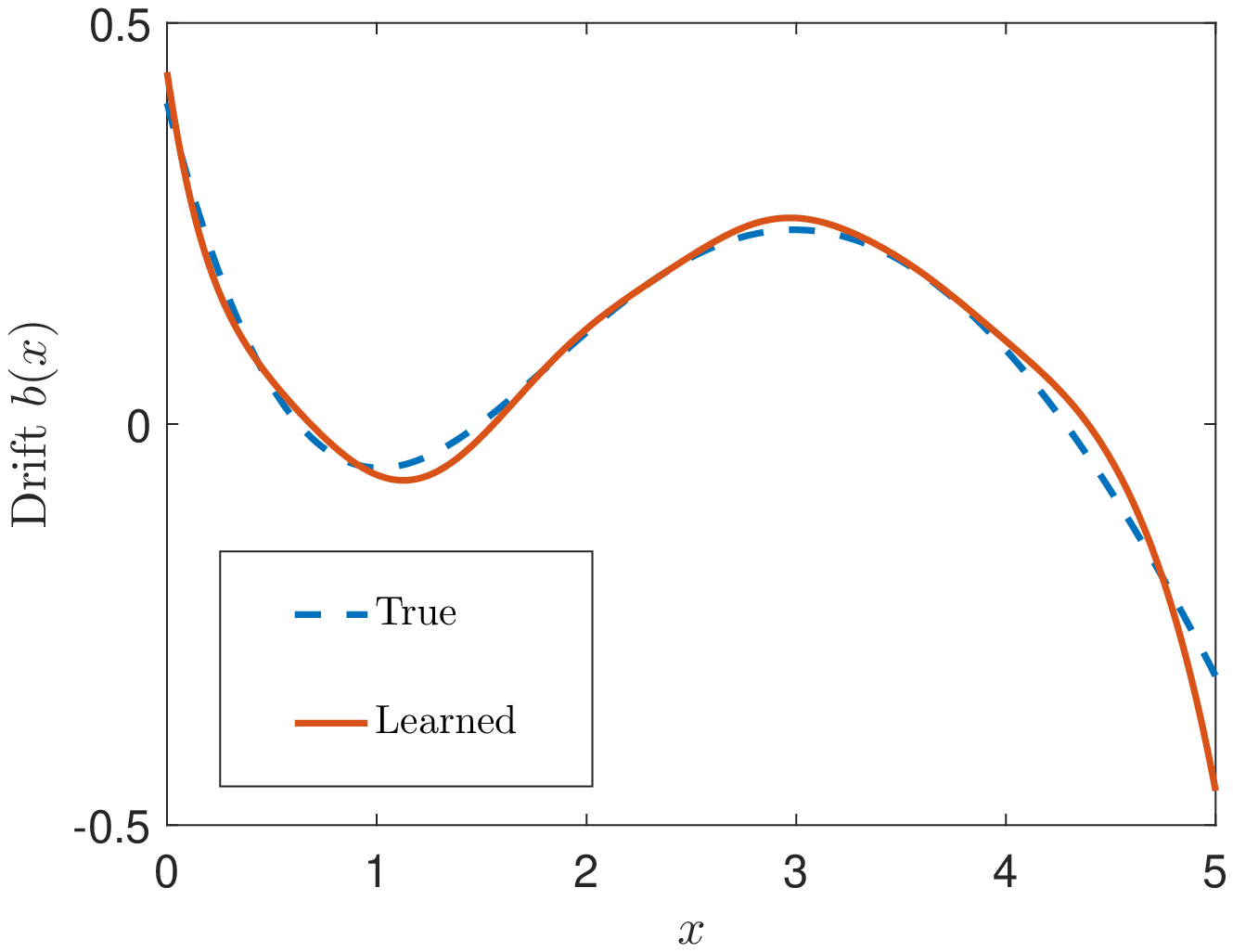}}
	\end{minipage}
	\hfill
	\begin{minipage}{0.48\linewidth}
		\leftline{(b)}
		\centerline{\includegraphics[height = 5cm, width = 6.5cm]{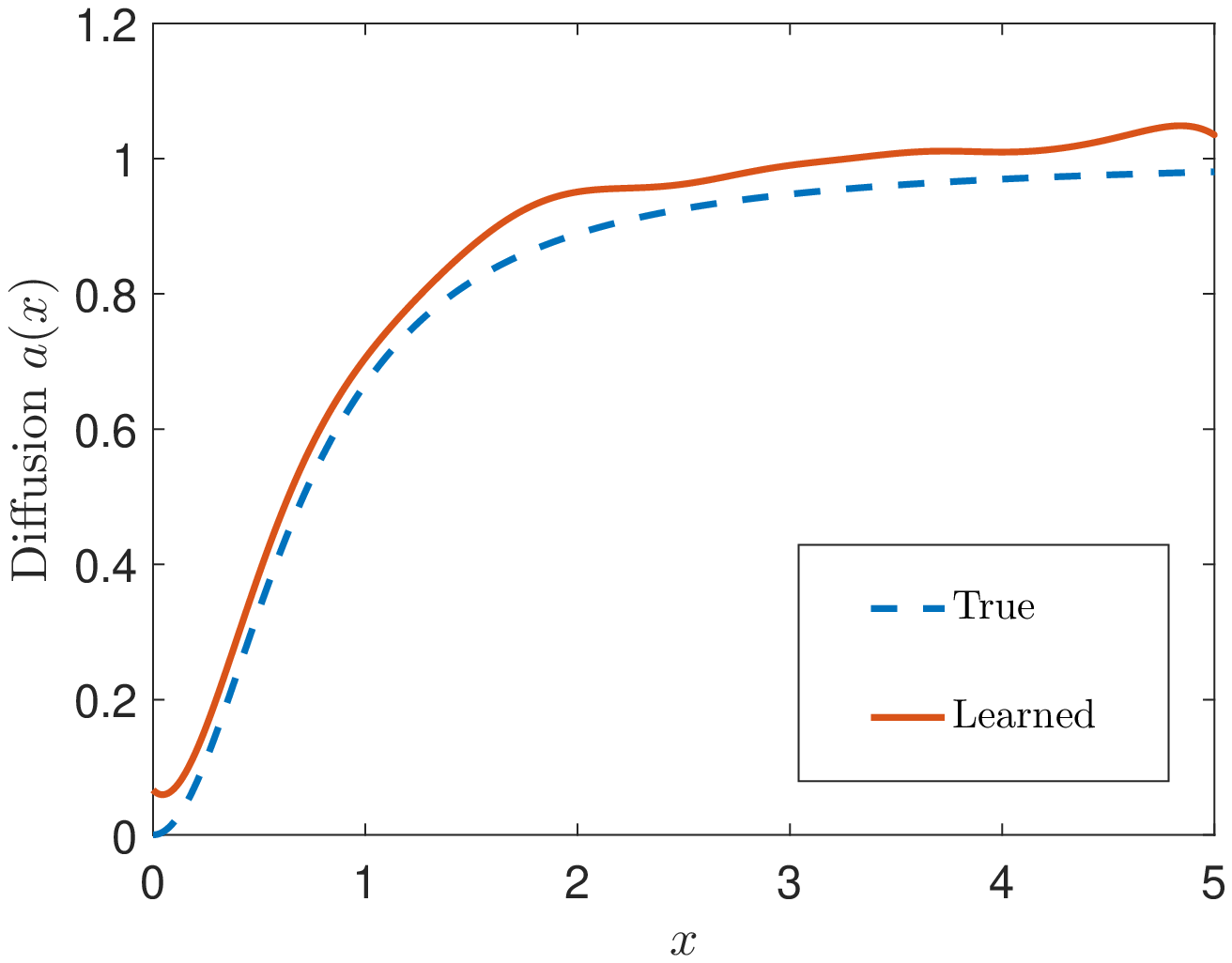}}
	\end{minipage}
	\caption{(a) True and learned drift coefficient. (b) True and learned diffusion coefficient.}\label{fig:3}
\end{figure}
\end{exa}

\section{Discussion}
\label{sec:5}
In this paper, we have designed a machine learning method to discover stochastic dynamical systems with both non-Gaussian L\'evy noise and Gaussian noise (Brownian motion), from observational, experimental or simulated data sets. We have generalized our recently proposed method from the special stochastic dynamical systems with  symmetric  L\'evy motion to the general systems with  asymmetric L\'evy motion.  Based on the expressions (i.e., nonlocal Kramers-Moyal formulas) of the jump measure, drift and diffusion coefficient in terms of  sample paths, we have further devised  a numerical algorithm for our machine learning method. The verification of this method on prototypical systems shows its efficacy and accuracy. Our novel approach provides a data-driven tool to extract stochastic governing laws for complex phenomena,  under general  non-Gaussian fluctuations.

Comparing with our previously published method, there are three  differences in this present work. First, it is generalized from   stochastic systems with symmetric  L\'evy motion to systems with  asymmetric  L\'evy motion.  The symmetric case is only an idealized situation. Thus what we should identify are not only the stability parameter and the L\'evy noise intensity but also the skewness parameter. The   application range of our numerical algorithm is thus much wider. Second, through removing the dependence of the different L\'evy components, the first assertion of Theorem \ref{thm1} is transformed into one-dimensional computation via the marginal probability distribution. Finally, the integration domain in the second and third assertions of Theorem \ref{thm1} is altered from a ball to a cube due to the independence of the  different L\'evy  components.

Remark that it is still a challenge to extend our approach to identify the stochastic governing equations with multiplicative L\'evy noise. In this case, the multiplicative L\'evy noise destroys the “space homogeneity” in the first assertion of Theorem \ref{thm1}, leading to the situation that the function $W$ depends on both $\mathbf{x}$ and $\mathbf{z}$ individually (not just depends only on the spatial dislocation $\mathbf{x}-\mathbf{z}$).

\section*{Acknowledgement}
This research was supported by the National Natural Science Foundation of China (Grants 11771449).

\begin{appendices}

 \section{$\alpha$-stable L\'evy processes}
\label{App:A}
A scalar $\alpha$-stable L\'evy process ${{L}_{t}}$ is a stochastic process with the following conditions:
\begin{enumerate}[i)]
\item ${{L}_{0}}=0$, a.s.;
\item Independent increments: for any choice of $n\ge 1$ and ${{t}_{0}}<{{t}_{1}}<\cdots <{{t}_{n-1}}<{{t}_{n}}$, the random variables ${{L}_{{{t}_{0}}}}$, ${{L}_{{{t}_{1}}}}-{{L}_{{{t}_{0}}}}$, ${{L}_{{{t}_{2}}}}-{{L}_{{{t}_{1}}}}$, $\cdots $, ${{L}_{{{t}_{n}}}}-{{L}_{{{t}_{n-1}}}}$ are independent;
\item Stationary increments: ${{L}_{t}}-{{L}_{s}}\sim {{S}_{\alpha }}\left( {{\left( t-s \right)}^{{1}/{\alpha }\;}},\beta ,0 \right)$;
\item Stochastically continuous sample paths: for every $s>0$, ${{L}_{t}}\to {{L}_{s}}$ in probability, as $t\to s$.
\end{enumerate}

The $\alpha$-stable L\'evy motion is a special but most popular type of the L\'evy process defined by the stable random variable with the distribution ${{S}_{\alpha }}\left( \delta ,\ \beta ,\ \lambda  \right)$ \cite{Ref31, Ref32, Ref33}. Usually, $\alpha \in \left( 0,\ 2 \right]$ is called the stability parameter, $\delta \in \left( 0,\ \infty  \right)$ is the scaling parameter, $\beta \in \left[ -1,\ 1 \right]$ is the skewness parameter and $\lambda \in \left( -\infty ,\ \infty  \right)$ is the shift parameter.

A stable random variable $X$ with $0<\alpha <2$ has the following “heavy tail” estimate:
$$
\underset{x\to \infty }{\mathop{\lim }}\,{{y}^{\alpha }}\mathbbm{P}\left( X>y \right)={{C}_{\alpha }}\frac{1+\beta }{2}{{\delta }^{\alpha }};
$$
where ${{C}_{\alpha }}$ is a positive constant depending on $\alpha$. In other words, the tail estimate decays polynomially. The $\alpha$-stable L\'evy motion has larger jumps with lower jump frequencies for smaller $\alpha$ ($0<\alpha <1$), while it has smaller jump sizes with higher jump frequencies for larger $\alpha$ ($1<\alpha <2$). The special case $\alpha=2$ corresponds to (Gaussian) Brownian motion. For more information aboutL\'evy process, refer to Refs. \cite{Ref29, Ref34}.

\section{Proof of Theorem 1}
\label{App:B}
i)  Since $p\left( \mathbf{x},0|\mathbf{z},0 \right)=\delta \left( \mathbf{x}-\mathbf{z} \right)$, then ${{p}_{i}}\left( {{x}_{i}},t|{{z}_{i}},0 \right)=\delta \left( {{x}_{i}}-{{z}_{i}} \right)$ and ${{p}_{i}}\left( {{x}_{i}},0|{{z}_{i}},0 \right)=0$ for arbitrary ${{x}_{i}}$ and  ${{z}_{i}}$ satisfying $\left| {{x}_{i}}-{{z}_{i}} \right|>\varepsilon $. For the convenience of representations, denote ${{\mathbf{x}}^{i}}$ as the rest of the vector $\mathbf{x}$ with  ${{x}_{i}}$ being removed. Then
$$
\begin{aligned}
  & \underset{t\to 0}{\mathop{\lim }}\,{{t}^{-1}}{{p}_{i}}\left( {{x}_{i}},t|{{z}_{i}},0 \right) \\
 & =\underset{t\to 0}{\mathop{\lim }}\,{{t}^{-1}}\left[ {{p}_{i}}\left( {{x}_{i}},t|{{z}_{i}},0 \right)-{{p}_{i}}\left( {{x}_{i}},0|{{z}_{i}},0 \right) \right] \\
 & =\int_{{{\mathbbm{R}}^{n-1}}}{{{\left. \frac{\partial p\left( \mathbf{x},t|\mathbf{z},0 \right)}{\partial t} \right|}_{t=0}}\textrm{d}{{\mathbf{x}}^{i}}}\\
 & =-\int_{{{\mathbbm{R}}^{n-1}}}{\nabla \cdot \left[ \mathbf{b}p\left( \mathbf{x},0|\mathbf{z},0 \right) \right]\textrm{d}{{\mathbf{x}}^{i}}}+\frac{1}{2}\int_{{{\mathbbm{R}}^{n-1}}}{Tr\left[ H\left( ap\left( \mathbf{x},0|\mathbf{z},0 \right) \right) \right]\textrm{d}{{\mathbf{x}}^{i}}} \\
 & \ \ -\sum\limits_{j=1}^{n}{\int_{{{\mathbbm{R}}^{n-1}}}{\int_{\mathbbm{R}\backslash \left\{ 0 \right\}}{\left[ p\left( \mathbf{x},0|\mathbf{z},0 \right)-p\left( \mathbf{x}-{{\sigma }_{j}}{y}_{j}{{\mathbf{e}}_{j}},0|\mathbf{z},0 \right)-{{\sigma }_{j}}\chi _{j}^{\alpha }\left( {y}_{j} \right){y}_{j}\frac{\partial }{\partial {{x}_{j}}}p\left( \mathbf{x},0|\mathbf{z},0 \right) \right]W_{j}^{\alpha ,\beta }\left( {y}_{j} \right)\textrm{d}{y}_{j}}\textrm{d}{{\mathbf{x}}^{i}}}}.
\end{aligned}
$$
Since $\left| {{x}_{i}}-{{z}_{i}} \right|>\varepsilon $, $p\left( \mathbf{x},0|\mathbf{z},0 \right)=0$. Thus the above equation is reduced as
$$
\sum\limits_{j=1}^{n}{\int_{{{\mathbbm{R}}^{n-1}}}{\int_{\mathbbm{R}\backslash \left\{ 0 \right\}}{p\left( \mathbf{x}-{{\sigma }_{j}}{y}_{j}{{\mathbf{e}}_{j}},0|\mathbf{z},0 \right)W_{j}^{\alpha ,\beta }\left( {y}_{j} \right)\textrm{d}{y}_{j}}\textrm{d}{{\mathbf{x}}^{i}}}}.
$$
Since
$$
\begin{aligned}
  & p\left( \mathbf{x}-{{\sigma }_{j}}{y}_{j}{{\mathbf{e}}_{j}},0|\mathbf{z},0 \right) \\
 & =\sigma _{j}^{-1}\delta \left( {{x}_{1}}-{{z}_{1}} \right)\cdots \delta \left( {{x}_{j-1}}-{{z}_{j-1}} \right)\delta \left( {{x}_{j}}-{{\sigma }_{j}}{y}_{j}-{{z}_{j}} \right)\delta \left( {{x}_{j+1}}-{{z}_{j+1}} \right)\cdots \delta \left( {{x}_{n}}-{{z}_{n}} \right), \\
\end{aligned}
$$
its integral is equal to zero for $j\ne i$ due to $\delta \left( {{x}_{i}}-{{z}_{i}} \right)=0$. Hence, we have
$$
\sum\limits_{j=1}^{n}{\int_{{{\mathbbm{R}}^{n-1}}}{\int_{\mathbbm{R}\backslash \left\{ 0 \right\}}{p\left( \mathbf{x}-{{\sigma }_{j}}{y}_{j}{{\mathbf{e}}_{j}},0|\mathbf{z},0 \right)W_{j}^{\alpha ,\beta }\left( {y}_{j} \right)\textrm{d}{y}_{j}}\textrm{d}{{\mathbf{x}}^{i}}}}=\sigma _{i}^{-1}W_{i}^{\alpha ,\beta }\left( \sigma _{i}^{-1}\left( {{x}_{i}}-{{z}_{i}} \right) \right).
$$

\noindent ii)  According to the Fokker-Planck equation (\ref{eq2}), we have
$$
\begin{aligned}
  & \underset{t\to 0}{\mathop{\lim }}\,{{t}^{-1}}\int_{\mathbf{x}-\mathbf{z}\in \Gamma }{\left( {{x}_{i}}-{{z}_{i}} \right)p\left( \mathbf{x},t|\mathbf{z},0 \right)\textrm{d}\mathbf{x}} \\
 & =\underset{t\to 0}{\mathop{\lim }}\,{{t}^{-1}}\int_{\mathbf{x}-\mathbf{z}\in \Gamma }{\left( {{x}_{i}}-{{z}_{i}} \right)\left[ p\left( \mathbf{x},t|\mathbf{z},0 \right)-p\left( \mathbf{x},0|\mathbf{z},0 \right)+p\left( \mathbf{x},0|\mathbf{z},0 \right) \right]\textrm{d}\mathbf{x}} \\
 & =\int_{\mathbf{x}-\mathbf{z}\in \Gamma }{\left( {{x}_{i}}-{{z}_{i}} \right){{\left. \frac{\partial p\left( \mathbf{x},t|\mathbf{z},0 \right)}{\partial t} \right|}_{t=0}}\textrm{d}\mathbf{x}}+\underset{t\to 0}{\mathop{\lim }}\,\left[ {{t}^{-1}}\int_{\mathbf{x}-\mathbf{z}\in \Gamma }{\left( {{x}_{i}}-{{z}_{i}} \right)\delta \left( \mathbf{x}-\mathbf{z} \right)\textrm{d}\mathbf{x}} \right] \\
 & =-\sum\limits_{j=1}^{n}{\int_{\mathbf{x}-\mathbf{z}\in \Gamma }{\left( {{x}_{i}}-{{z}_{i}} \right)\frac{\partial }{\partial {{x}_{j}}}\left[ {{b}_{j}}\left( \mathbf{x} \right)p\left( \mathbf{x},0|\mathbf{z},0 \right) \right]\textrm{d}\mathbf{x}}} \\
 & \ \ \ +\frac{1}{2}\sum\limits_{k,l=1}^{n}{\int_{\mathbf{x}-\mathbf{z}\in \Gamma }{\left( {{x}_{i}}-{{z}_{i}} \right)\frac{{{\partial }^{2}}}{\partial {{x}_{k}}\partial {{x}_{l}}}\left[ {{a}_{kl}}\left( \mathbf{x} \right)p\left( \mathbf{x},0|\mathbf{z},0 \right) \right]\textrm{d}\mathbf{x}}} \\
 & \ \ \ -\sum\limits_{j=1}^{n}{\int_{\mathbf{x}-\mathbf{z}\in \Gamma }{\int_{\mathbbm{R}\backslash \left\{ 0 \right\}}{\left( {{x}_{i}}-{{z}_{i}} \right)\left[ p\left( \mathbf{x},0|\mathbf{z},0 \right)-p\left( \mathbf{x}-{{\sigma }_{j}}{y}_{j}{{\mathbf{e}}_{j}},0|\mathbf{z},0 \right)-{{\sigma }_{j}}\chi _{j}^{\alpha }\left( {y}_{j} \right){y}_{j}\frac{\partial }{\partial {{x}_{j}}}p\left( \mathbf{x},0|\mathbf{z},0 \right) \right]W_{j}^{\alpha ,\beta }\left( {y}_{j} \right)\textrm{d}{y}_{j}}\textrm{d}\mathbf{x}}}. \\
\end{aligned}
$$
The application of integration by parts into the first term leads to
$$
\begin{aligned}
  & \sum\limits_{j=1}^{n}{\int_{\mathbf{x}-\mathbf{z}\in \Gamma }{\left( {{x}_{i}}-{{z}_{i}} \right)\frac{\partial }{\partial {{x}_{j}}}\left[ {{b}_{j}}\left( \mathbf{x} \right)p\left( \mathbf{x},0|\mathbf{z},0 \right) \right]\textrm{d}\mathbf{x}}} \\
 & =-\sum\limits_{j=1}^{n}{\int_{\mathbf{x}-\mathbf{z}\in \Gamma }{{{b}_{j}}\left( \mathbf{x} \right)p\left( \mathbf{x},0|\mathbf{z},0 \right)\frac{\partial }{\partial {{x}_{j}}}\left( {{x}_{i}}-{{z}_{i}} \right)\textrm{d}\mathbf{x}}} \\
 & =-\sum\limits_{j=1}^{n}{\int_{\mathbf{x}-\mathbf{z}\in \Gamma }{{{b}_{j}}\left( \mathbf{x} \right)\delta \left( \mathbf{x}-\mathbf{z} \right){{\delta }_{ij}}\textrm{d}\mathbf{x}}} \\
 & =-{{b}_{i}}\left( \mathbf{z} \right). \\
\end{aligned}
$$
Therein, the boundary condition vanishes since $p\left( \mathbf{x},0|\mathbf{z},0 \right)=0$ as $\left| {{x}_{j}}-{{z}_{j}} \right|=\varepsilon $. For the second integration, we use integration by parts twice
$$
\begin{aligned}
  & \sum\limits_{k,l=1}^{n}{\int_{\mathbf{x}-\mathbf{z}\in \Gamma }{\left( {{x}_{i}}-{{z}_{i}} \right)\frac{{{\partial }^{2}}}{\partial {{x}_{k}}\partial {{x}_{l}}}\left[ {{a}_{kl}}\left( \mathbf{x} \right)p\left( \mathbf{x},0|\mathbf{z},0 \right) \right]\textrm{d}\mathbf{x}}} \\
 & =-\sum\limits_{k,l=1}^{n}{\int_{\mathbf{x}-\mathbf{z}\in \Gamma }{\frac{\partial }{\partial {{x}_{k}}}\left( {{x}_{i}}-{{z}_{i}} \right)\frac{\partial }{\partial {{x}_{l}}}\left[ {{a}_{kl}}\left( \mathbf{x} \right)p\left( \mathbf{x},0|\mathbf{z},0 \right) \right]\textrm{d}\mathbf{x}}} \\
 & =\sum\limits_{k,l=1}^{n}{\int_{\mathbf{x}-\mathbf{z}\in \Gamma }{\left( {{x}_{i}}-{{z}_{i}} \right){{a}_{kl}}\left( \mathbf{x} \right)p\left( \mathbf{x},0|\mathbf{z},0 \right)\frac{\partial {{\delta }_{ik}}}{\partial {{x}_{l}}}\textrm{d}\mathbf{x}}} \\
 & =0. \\
\end{aligned}
$$
For the third integration, we derive it separately. First, according to Tonelli’s theorem \cite{Ref35}, we obtain
$$
\begin{aligned}
  & \sum\limits_{j=1}^{n}{\int_{\mathbf{x}-\mathbf{z}\in \Gamma }{\int_{\mathbbm{R}\backslash \left\{ 0 \right\}}{\left( {{x}_{i}}-{{z}_{i}} \right)p\left( \mathbf{x},0|\mathbf{z},0 \right)W_{j}^{\alpha ,\beta }\left( {y}_{j} \right)\textrm{d}{y}_{j}}\textrm{d}\mathbf{x}}} \\
 & =\sum\limits_{j=1}^{n}{\int_{\mathbbm{R}\backslash \left\{ 0 \right\}}{\int_{\mathbf{x}-\mathbf{z}\in \Gamma }{\left( {{x}_{i}}-{{z}_{i}} \right)\delta \left( \mathbf{x}-\mathbf{z} \right)\textrm{d}\mathbf{x}}W_{j}^{\alpha ,\beta }\left( {y}_{j} \right)\textrm{d}{y}_{j}}} \\
 & =0. \\
\end{aligned}
$$
Second,
$$
\begin{aligned}
  & \sum\limits_{j=1}^{n}{\int_{\mathbf{x}-\mathbf{z}\in \Gamma }{\int_{\mathbbm{R}\backslash \left\{ 0 \right\}}{\left( {{x}_{i}}-{{z}_{i}} \right)p\left( \mathbf{x}-{{\sigma }_{j}}{y}_{j}{{\mathbf{e}}_{j}},0|\mathbf{z},0 \right)W_{j}^{\alpha ,\beta }\left( {y}_{j} \right)\textrm{d}{y}_{j}}\textrm{d}\mathbf{x}}} \\
 & =\int_{\mathbf{x}-\mathbf{z}\in \Gamma }{\int_{\mathbbm{R}\backslash \left\{ 0 \right\}}{\left( {{x}_{i}}-{{z}_{i}} \right)p\left( \mathbf{x}-{{\sigma }_{i}}{y}_{i}{{\mathbf{e}}_{i}},0|\mathbf{z},0 \right)W_{i}^{\alpha ,\beta }\left( {y}_{i} \right)\textrm{d}{y}_{i}}\textrm{d}\mathbf{x}} \\
 & =\sigma _{i}^{-1}\int_{-\varepsilon }^{\varepsilon }{{y}_{i}W_{i}^{\alpha ,\beta }\left( \sigma _{i}^{-1}{y}_{i} \right)\textrm{d}{y}_{i}}. \\
\end{aligned}
$$
Finally, according to Tonelli’s theorem and integration by parts, we have
$$
\begin{aligned}
  & \sum\limits_{j=1}^{n}{\int_{\mathbf{x}-\mathbf{z}\in \Gamma }{\int_{\mathbbm{R}\backslash \left\{ 0 \right\}}{\left( {{x}_{i}}-{{z}_{i}} \right){{\sigma }_{j}}\chi _{j}^{\alpha }\left( {y}_{j} \right){y}_{j}\frac{\partial }{\partial {{x}_{j}}}p\left( \mathbf{x},0|\mathbf{z},0 \right)W_{j}^{\alpha ,\beta }\left( {y}_{j} \right)\textrm{d}{y}_{j}}\textrm{d}\mathbf{x}}} \\
 & =\sum\limits_{j=1}^{n}{\int_{\mathbbm{R}\backslash \left\{ 0 \right\}}{\int_{\mathbf{x}-\mathbf{z}\in \Gamma }{\left( {{x}_{i}}-{{z}_{i}} \right)\frac{\partial }{\partial {{x}_{j}}}p\left( \mathbf{x},0|\mathbf{z},0 \right)\textrm{d}\mathbf{x}}{{\sigma }_{j}}\chi _{j}^{\alpha }\left( {y}_{j} \right){y}_{j}W_{j}^{\alpha ,\beta }\left( {y}_{j} \right)\textrm{d}{y}_{j}}} \\
 & =-\sum\limits_{j=1}^{n}{\int_{\mathbbm{R}\backslash \left\{ 0 \right\}}{\int_{\mathbf{x}-\mathbf{z}\in \Gamma }{p\left( \mathbf{x},0|\mathbf{z},0 \right)\frac{\partial }{\partial {{x}_{j}}}\left( {{x}_{i}}-{{z}_{i}} \right)\textrm{d}\mathbf{x}}{{\sigma }_{j}}\chi _{j}^{\alpha }\left( {y}_{j} \right){y}_{j}W_{j}^{\alpha ,\beta }\left( {y}_{j} \right)\textrm{d}{y}_{j}}} \\
 & =\left\{ \begin{array}{ccc}
   0, & \alpha <1,  \\
   -\sigma _{i}^{-1}\int_{-1}^{1}{{y}_{i}W_{i}^{\alpha ,\beta }\left( \sigma _{i}^{-1}{y}_{i} \right)\textrm{d}{y}_{i}}, & \alpha =1,  \\
   -\sigma _{i}^{-1}\int_{-\infty }^{\infty }{{y}_{i}W_{i}^{\alpha ,\beta }\left( \sigma _{i}^{-1}{y}_{i} \right)\textrm{d}{y}_{i}}, & \alpha >1  \\
\end{array} \right.. \\
\end{aligned}
$$
Hence,
$$
\underset{t\to 0}{\mathop{\lim }}\,{{t}^{-1}}\int_{\mathbf{x}-\mathbf{z}\in \Gamma }{\left( {{x}_{i}}-{{z}_{i}} \right)p\left( \mathbf{x},t|\mathbf{z},0 \right)\textrm{d}\mathbf{x}}={{b}_{i}}\left( \mathbf{z} \right)+R_{i}^{\alpha ,\beta }\left( \varepsilon  \right).
$$

\noindent iii) According to the Fokker-Planck equation (\ref{eq2}), we have
$$
\begin{aligned}
  & \underset{t\to 0}{\mathop{\lim }}\,{{t}^{-1}}\int_{\mathbf{x}-\mathbf{z}\in \Gamma }{\left( {{x}_{i}}-{{z}_{i}} \right)\left( {{x}_{j}}-{{z}_{j}} \right)p\left( \mathbf{x},t|\mathbf{z},0 \right)\textrm{d}\mathbf{x}} \\
 & =\underset{t\to 0}{\mathop{\lim }}\,{{t}^{-1}}\int_{\mathbf{x}-\mathbf{z}\in \Gamma }{\left( {{x}_{i}}-{{z}_{i}} \right)\left( {{x}_{j}}-{{z}_{j}} \right)\left[ p\left( \mathbf{x},t|\mathbf{z},0 \right)-p\left( \mathbf{x},0|\mathbf{z},0 \right)+p\left( \mathbf{x},0|\mathbf{z},0 \right) \right]\textrm{d}\mathbf{x}} \\
 & =\int_{\mathbf{x}-\mathbf{z}\in \Gamma }{\left( {{x}_{i}}-{{z}_{i}} \right)\left( {{x}_{j}}-{{z}_{j}} \right){{\left. \frac{\partial p\left( \mathbf{x},t|\mathbf{z},0 \right)}{\partial t} \right|}_{t=0}}\textrm{d}\mathbf{x}} \\
 & \ \ \ +\underset{t\to 0}{\mathop{\lim }}\,\left[ {{t}^{-1}}\int_{\mathbf{x}-\mathbf{z}\in \Gamma }{\left( {{x}_{i}}-{{z}_{i}} \right)\left( {{x}_{j}}-{{z}_{j}} \right)\delta \left( \mathbf{x}-\mathbf{z} \right)\textrm{d}\mathbf{x}} \right] \\
 & =-\sum\limits_{k=1}^{n}{\int_{\mathbf{x}-\mathbf{z}\in \Gamma }{\left( {{x}_{i}}-{{z}_{i}} \right)\left( {{x}_{j}}-{{z}_{j}} \right)\frac{\partial }{\partial {{x}_{k}}}\left[ {{b}_{k}}\left( \mathbf{x} \right)p\left( \mathbf{x},0|\mathbf{z},0 \right) \right]\textrm{d}\mathbf{x}}} \\
 & \ \ \ +\frac{1}{2}\sum\limits_{k,l=1}^{n}{\int_{\mathbf{x}-\mathbf{z}\in \Gamma }{\left( {{x}_{i}}-{{z}_{i}} \right)\left( {{x}_{j}}-{{z}_{j}} \right)\frac{{{\partial }^{2}}}{\partial {{x}_{k}}\partial {{x}_{l}}}\left[ {{a}_{kl}}\left( \mathbf{x} \right)p\left( \mathbf{x},0|\mathbf{z},0 \right) \right]\textrm{d}\mathbf{x}}} \\
 & \ \ \ -\sum\limits_{k=1}^{n}{\int_{\mathbf{x}-\mathbf{z}\in \Gamma }{\int_{\mathbbm{R}\backslash \left\{ 0 \right\}}{\left( {{x}_{i}}-{{z}_{i}} \right)\left( {{x}_{j}}-{{z}_{j}} \right)\left[ p\left( \mathbf{x},0|\mathbf{z},0 \right)-p\left( \mathbf{x}-{{\sigma }_{k}}{y}_{k}{{\mathbf{e}}_{k}},0|\mathbf{z},0 \right)-{{\sigma }_{k}}\chi _{k}^{\alpha }\left( {y}_{k} \right){y}_{k}\frac{\partial }{\partial {{x}_{k}}}p\left( \mathbf{x},0|\mathbf{z},0 \right) \right]W_{k}^{\alpha ,\beta }\left( {y}_{k} \right)\textrm{d}{y}_{k}}\textrm{d}\mathbf{x}}}. \\
\end{aligned}
$$
The application of integration by parts into the first term yields
$$
\begin{aligned}
  & \sum\limits_{k=1}^{n}{\int_{\mathbf{x}-\mathbf{z}\in \Gamma }{\left( {{x}_{i}}-{{z}_{i}} \right)\left( {{x}_{j}}-{{z}_{j}} \right)\frac{\partial }{\partial {{x}_{k}}}\left[ {{b}_{k}}\left( \mathbf{x} \right)p\left( \mathbf{x},0|\mathbf{z},0 \right) \right]\textrm{d}\mathbf{x}}} \\
 & =-\sum\limits_{k=1}^{n}{\int_{\mathbf{x}-\mathbf{z}\in \Gamma }{{{b}_{k}}\left( \mathbf{x} \right)p\left( \mathbf{x},0|\mathbf{z},0 \right)\frac{\partial }{\partial {{x}_{k}}}\left[ \left( {{x}_{i}}-{{z}_{i}} \right)\left( {{x}_{j}}-{{z}_{j}} \right) \right]\textrm{d}\mathbf{x}}} \\
 & =-\sum\limits_{k=1}^{n}{\int_{\mathbf{x}-\mathbf{z}\in \Gamma }{{{b}_{k}}\left( \mathbf{x} \right)\delta \left( \mathbf{x}-\mathbf{z} \right)\left[ {{\delta }_{ik}}\left( {{x}_{j}}-{{z}_{j}} \right)+{{\delta }_{jk}}\left( {{x}_{i}}-{{z}_{i}} \right) \right]\textrm{d}\mathbf{x}}} \\
 & =0. \\
\end{aligned}
$$
Therein, the boundary condition vanishes since $p\left( \mathbf{x},0|\mathbf{z},0 \right)=0$ as $\left| {{x}_{k}}-{{z}_{k}} \right|=\varepsilon $. For the second integration, we use integration by parts again
$$
\begin{aligned}
  & \frac{1}{2}\sum\limits_{k,l=1}^{n}{\int_{\mathbf{x}-\mathbf{z}\in \Gamma }{\left( {{x}_{i}}-{{z}_{i}} \right)\left( {{x}_{j}}-{{z}_{j}} \right)\frac{{{\partial }^{2}}}{\partial {{x}_{k}}\partial {{x}_{l}}}\left[ {{a}_{kl}}\left( \mathbf{x} \right)p\left( \mathbf{x},0|\mathbf{z},0 \right) \right]\textrm{d}\mathbf{x}}} \\
 & =-\frac{1}{2}\sum\limits_{k,l=1}^{n}{\int_{\mathbf{x}-\mathbf{z}\in \Gamma }{\frac{\partial }{\partial {{x}_{k}}}\left[ \left( {{x}_{i}}-{{z}_{i}} \right)\left( {{x}_{j}}-{{z}_{j}} \right) \right]\frac{\partial }{\partial {{x}_{l}}}\left[ {{a}_{kl}}\left( \mathbf{x} \right)p\left( \mathbf{x},0|\mathbf{z},0 \right) \right]\textrm{d}\mathbf{x}}} \\
 & =\frac{1}{2}\sum\limits_{k,l=1}^{n}{\int_{\mathbf{x}-\mathbf{z}\in \Gamma }{{{a}_{kl}}\left( \mathbf{x} \right)p\left( \mathbf{x},0|\mathbf{z},0 \right)\frac{\partial }{\partial {{x}_{l}}}\left[ {{\delta }_{ik}}\left( {{x}_{j}}-{{z}_{j}} \right)+{{\delta }_{jk}}\left( {{x}_{i}}-{{z}_{i}} \right) \right]\textrm{d}\mathbf{x}}} \\
 & =\frac{1}{2}\sum\limits_{k,l=1}^{n}{\int_{\mathbf{x}-\mathbf{z}\in \Gamma }{{{a}_{kl}}\left( \mathbf{x} \right)\delta \left( \mathbf{x}-\mathbf{z} \right)\left( {{\delta }_{ik}}{{\delta }_{jl}}+{{\delta }_{il}}{{\delta }_{jk}} \right)\textrm{d}\mathbf{x}}} \\
 & =\frac{1}{2}\left[ {{a}_{ij}}\left( \mathbf{z} \right)+{{a}_{ji}}\left( \mathbf{z} \right) \right] \\
 & ={{a}_{ij}}\left( \mathbf{z} \right). \\
\end{aligned}
$$
We still derive the third integration separately. First, according to Tonelli’s theorem, we obtain
$$
\begin{aligned}
  & \sum\limits_{j=1}^{n}{\int_{\mathbf{x}-\mathbf{z}\in \Gamma }{\int_{\mathbbm{R}\backslash \left\{ 0 \right\}}{\left( {{x}_{i}}-{{z}_{i}} \right)\left( {{x}_{j}}-{{z}_{j}} \right)p\left( \mathbf{x},0|\mathbf{z},0 \right)W_{j}^{\alpha ,\beta }\left( {y}_{j} \right)\textrm{d}{y}_{j}}\textrm{d}\mathbf{x}}} \\
 & =\sum\limits_{j=1}^{n}{\int_{\mathbbm{R}\backslash \left\{ 0 \right\}}{\int_{\mathbf{x}-\mathbf{z}\in \Gamma }{\left( {{x}_{i}}-{{z}_{i}} \right)\left( {{x}_{j}}-{{z}_{j}} \right)\delta \left( \mathbf{x}-\mathbf{z} \right)\textrm{d}\mathbf{x}}W_{j}^{\alpha ,\beta }\left( {y}_{j} \right)\textrm{d}{y}_{j}}} \\
 & =0. \\
\end{aligned}
$$
Second, for $i\ne j$,
$$
\sum\limits_{k=1}^{n}{\int_{\mathbf{x}-\mathbf{z}\in \Gamma }{\int_{\mathbbm{R}\backslash \left\{ 0 \right\}}{\left( {{x}_{i}}-{{z}_{i}} \right)\left( {{x}_{j}}-{{z}_{j}} \right)p\left( \mathbf{x}-{{\sigma }_{k}}{y}_{k}{{\mathbf{e}}_{k}},0|\mathbf{z},0 \right)W_{k}^{\alpha ,\beta }\left( {y}_{k} \right)\textrm{d}{y}_{k}}\textrm{d}\mathbf{x}}}=0
$$
and for $i=j$,
$$
\begin{aligned}
  & \sum\limits_{k=1}^{n}{\int_{\mathbf{x}-\mathbf{z}\in \Gamma }{\int_{\mathbbm{R}\backslash \left\{ 0 \right\}}{{{\left( {{x}_{i}}-{{z}_{i}} \right)}^{2}}p\left( \mathbf{x}-{{\sigma }_{k}}{y}_{k}{{\mathbf{e}}_{k}},0|\mathbf{z},0 \right)W_{k}^{\alpha ,\beta }\left( {y}_{k} \right)\textrm{d}{y}_{k}}\textrm{d}\mathbf{x}}} \\
 & =\int_{\mathbf{x}-\mathbf{z}\in \Gamma }{\int_{\mathbbm{R}\backslash \left\{ 0 \right\}}{{{\left( {{x}_{i}}-{{z}_{i}} \right)}^{2}}p\left( \mathbf{x}-{{\sigma }_{i}}{y}_{i}{{\mathbf{e}}_{i}},0|\mathbf{z},0 \right)W_{i}^{\alpha ,\beta }\left( {y}_{i} \right)\textrm{d}{y}_{i}}\textrm{d}\mathbf{x}} \\
 & =\sigma _{i}^{-1}\int_{-\varepsilon }^{\varepsilon }{{{y}_{i}^{2}}W_{i}^{\alpha ,\beta }\left( \sigma _{i}^{-1}{y}_{i} \right)\textrm{d}{y}_{i}}. \\
\end{aligned}
$$
This integration is bounded according to the definition of the jump measure. Finally, according to Tonelli’s theorem and integration by parts, we have
$$
\begin{aligned}
  & \sum\limits_{k=1}^{n}{\int_{\mathbf{x}-\mathbf{z}\in \Gamma }{\int_{\mathbbm{R}\backslash \left\{ 0 \right\}}{\left( {{x}_{i}}-{{z}_{i}} \right)\left( {{x}_{j}}-{{z}_{j}} \right){{\sigma }_{k}}\chi _{k}^{\alpha }{y}_{k}\frac{\partial }{\partial {{x}_{k}}}p\left( \mathbf{x},0|\mathbf{z},0 \right)W_{k}^{\alpha ,\beta }\left( {y}_{k} \right)\textrm{d}{y}_{k}}\textrm{d}\mathbf{x}}} \\
 & =\sum\limits_{k=1}^{n}{\int_{\mathbbm{R}\backslash \left\{ 0 \right\}}{\int_{\mathbf{x}-\mathbf{z}\in \Gamma }{\left( {{x}_{i}}-{{z}_{i}} \right)\left( {{x}_{j}}-{{z}_{j}} \right)\frac{\partial }{\partial {{x}_{k}}}p\left( \mathbf{x},0|\mathbf{z},0 \right)\textrm{d}\mathbf{x}}{{\sigma }_{k}}\chi _{k}^{\alpha }{y}_{k}W_{k}^{\alpha ,\beta }\left( {y}_{k}\right)\textrm{d}{y}_{k}}} \\
 & =-\sum\limits_{k=1}^{n}{\int_{\mathbbm{R}\backslash \left\{ 0 \right\}}{\int_{\mathbf{x}-\mathbf{z}\in \Gamma }{\delta \left( \mathbf{x}-\mathbf{z} \right)\left[ {{\delta }_{ik}}\left( {{x}_{j}}-{{z}_{j}} \right)+{{\delta }_{jk}}\left( {{x}_{i}}-{{z}_{i}} \right) \right]\textrm{d}\mathbf{x}}{{\sigma }_{k}}\chi _{k}^{\alpha }{y}_{k}W_{k}^{\alpha ,\beta }\left( {y}_{k} \right)\textrm{d}{y}_{k}}} \\
 & =0. \\
\end{aligned}
$$
Hence,
$$
\underset{t\to 0}{\mathop{\lim }}\,{{t}^{-1}}\int_{\mathbf{x}-\mathbf{z}\in \Gamma }{\left( {{x}_{i}}-{{z}_{i}} \right)\left( {{x}_{j}}-{{z}_{j}} \right)p\left( \mathbf{x},t|\mathbf{z},0 \right)\textrm{d}\mathbf{x}}={{a}_{ij}}\left( \mathbf{z} \right)+S_{ij}^{\alpha ,\beta }\left( \varepsilon  \right).
$$

\noindent The proof is complete.

\section{Proof of Corollary 2}
\label{App:C}
i)  This is derived directly by integrating the equation in the first assertion of Theorem \ref{thm1} about ${{x}_{i}}$ on the interval $\left[ {{c}_{1}},\ {{c}_{2}} \right)$.

\noindent ii)  Let the set $\textrm{d}U=\left[ {{u}_{1}},\ {{u}_{1}}+\textrm{d}{{u}_{1}} \right)\times \left[ {{u}_{2}},\ {{u}_{2}}+\textrm{d}{{u}_{2}} \right)\times \cdots \times \left[ {{u}_{n}},\ {{u}_{n}}+\textrm{d}{{u}_{n}} \right)$. Then we have
$$
\begin{aligned}
  & \mathbbm{P}\left\{ \left. \mathbf{x}\left( t \right)\in \textrm{d}U;\ \mathbf{x}\left( t \right)-\mathbf{z}\in \Gamma  \right|\mathbf{x}\left( 0 \right)=\mathbf{z} \right\} \\
 & =\mathbbm{P}\left\{ \left. \mathbf{x}\left( t \right)\in \textrm{d}U \right|\mathbf{x}\left( 0 \right)=\mathbf{z};\ \mathbf{x}\left( t \right)-\mathbf{z}\in \Gamma  \right\}\cdot \mathbbm{P}\left\{ \left. \mathbf{x}\left( t \right)-\mathbf{z}\in \Gamma  \right|\mathbf{x}\left( 0 \right)=\mathbf{z} \right\}. \\
\end{aligned}
$$
Thus
$$
\begin{aligned}
  & \int_{\mathbf{x}-\mathbf{z}\in \Gamma }{\left( {{x}_{i}}-{{z}_{i}} \right)p\left( \mathbf{x},t|\mathbf{z},0 \right)\textrm{d}\mathbf{x}} \\
 & =\int_{\mathbf{u}-\mathbf{z}\in \Gamma }{\left( {{u}_{i}}-{{z}_{i}} \right)\mathbbm{P}\left\{ \left. \mathbf{x}\left( t \right)\in \textrm{d}U \right|\mathbf{x}\left( 0 \right)=\mathbf{z} \right\}} \\
 & =\int_{\mathbf{u}-\mathbf{z}\in \Gamma }{\left( {{u}_{i}}-{{z}_{i}} \right)\mathbbm{P}\left\{ \left. \mathbf{x}\left( t \right)\in \textrm{d}U;\ \mathbf{x}\left( t \right)-\mathbf{z}\in \Gamma  \right|\mathbf{x}\left( 0 \right)=\mathbf{z} \right\}} \\
 & =\mathbbm{P}\left\{ \left. \mathbf{x}\left( t \right)-\mathbf{z}\in \Gamma  \right|\mathbf{x}\left( 0 \right)=\mathbf{z} \right\}\cdot \int_{\left| \mathbf{u}-\mathbf{z} \right|<\varepsilon }{\left( {{u}_{i}}-{{z}_{i}} \right)\mathbbm{P}\left\{ \left. \mathbf{x}\left( t \right)\in \textrm{d}U \right|\mathbf{x}\left( 0 \right)=\mathbf{z};\ \mathbf{x}\left( t \right)-\mathbf{z}\in \Gamma  \right\}} \\
 & =\mathbbm{P}\left\{ \left. \mathbf{x}\left( t \right)-\mathbf{z}\in \Gamma  \right|\mathbf{x}\left( 0 \right)=\mathbf{z} \right\}\cdot \textrm{E}\left[ \left. \left( {{x}_{i}}\left( t \right)-{{z}_{i}} \right) \right|\mathbf{x}\left( 0 \right)=\mathbf{z};\ \mathbf{x}\left( t \right)-\mathbf{z}\in \Gamma  \right]. \\
\end{aligned}
$$
Hence, the conclusion is immediately deduced
$$
\underset{t\to 0}{\mathop{\lim }}\,{{t}^{-1}}\mathbbm{P}\left\{ \left. \mathbf{x}\left( t \right)-\mathbf{z}\in \Gamma  \right|\mathbf{x}\left( 0 \right)=\mathbf{z} \right\}\cdot \textrm{E}\left[ \left. \left( {{x}_{i}}\left( t \right)-{{z}_{i}} \right) \right|\mathbf{x}\left( 0 \right)=\mathbf{z};\ \mathbf{x}\left( t \right)-\mathbf{z}\in \Gamma  \right]={{b}_{i}}\left( \mathbf{z} \right)+R_{i}^{\alpha ,\beta }\left( \varepsilon  \right).
$$

\noindent iii) This proof is similar to the second conclusion.

\noindent The proof is complete.

\end{appendices}

\section*{Data Availability Statement}
The data that support the findings of this study are openly available in GitHub \cite{Ref36}.


%
%


\begin{thebibliography}{}
%
%

\bibitem{Ref8}
L. Boninsegna, F. N\"uske, C. Clementi, Sparse learning of stochastic dynamical equations. J. Chem. Phys. \textbf{148}, 241723 (2018)

\bibitem{DaiMinChaos}
M. Dai, T. Gao, Y. Lu, Y. Zheng, J. Duan, Detecting the maximum likelihood transition path from data of stochastic dynamical systems. Chaos \textbf{30}, 113124 (2020)

\bibitem{KM}
M. Reza Rahimi Tabar, \emph{Analysis and Data-Based Reconstruction of Complex Nonlinear Dynamical Systems}, Springer (2019)

\bibitem{Ref1}
M. O. Williams, I. G. Kevrekidis, C. W. Rowley, A Data-Driven Approximation of the Koopman Operator: Extending Dynamic Mode Decomposition. J. Nonlinear Sci. \textbf{25}, 1307–1346 (2015)

\bibitem{Ref2}
S. Klus, F. Nüske, S. Peitz, J. H. Niemann, C. Clementi, C. Schütte, Data-driven approximation of the Koopman generator: Model reduction, system identification, and control. Phys. D Nonlinear Phenom. \textbf{406}, 132416 (2020)

\bibitem{WuFuDuan2019}
D. Wu, M. Fu, J. Duan, Discovering mean residence time and escape probability from data of stochastic dynamical systems. Chaos \textbf{29}, 093122 (2019)

\bibitem{LuYB2020}
Y. Lu, J. Duan, Discovering transition phenomena from data of stochastic dynamical systems with L\'evy noise. Chaos \textbf{30}, 093110 (2020)

\bibitem{Yanxia2020}
Y. Zhang,  J. Duan, Y. Jin, Y. Li, Extracting non-Gaussian governing laws from data on mean exit time. Chaos \textbf{30}, 113112 (2020)

\bibitem{Ref3}
S. L. Brunton, J. L. Proctor, J. N. Kutz, Discovering governing equations from data by sparse identification of nonlinear dynamical systems. Proc. Natl. Acad. Sci. \textbf{113}, 3932–3937 (2016)

\bibitem{Ref4}
R. González-García, R. Rico-Martínez, I. G. Kevrekidis, Identification of distributed parameter systems: A neural net based approach. Comput. Chem. Eng. \textbf{22}, S965–S968 (1998)

\bibitem{Ref5}
H. Schaeffer, R. Caflisch, C. D. Hauck, S. Osher, Sparse dynamics for partial differential equations. Proc. Natl. Acad. Sci. \textbf{110}, 6634–6639 (2013)

\bibitem{Ref6}
S. Rudy, A. Alla, S. L. B. S, J. N. Kutz, Data-Driven Identification of Parametric Partial Differential Equations. SIAM \textbf{18}, 643–660 (2019)

\bibitem{Ref7}
S. H. Rudy, S. L. Brunton, J. L. Proctor, J. N. Kutz, Data-driven discovery of partial differential equations. Sci. Adv. \textbf{3}, e1602614 (2017)

\bibitem{Barn}
O. E. Barndorff-Nielsen, T. Mikosch and S. I. Resnick (Eds.), \emph{L\'evy Processes: Theory and Applications}, Birkh\"auser, Boston (2001)

\bibitem{Ref9}
P. D. Ditlevsen, Observation of $\alpha$-stable noise induced millennial climate changes from an ice-core record. Geophys. Res. Lett. \textbf{26}, 1441–1444 (1999)

\bibitem{Ref10}
J. M. Raser, E. K. O. Shea, Noise in Gene Expression: Origins, Consequences, and Control. Science (80-. ). \textbf{309}, 2010–2013 (2005)

\bibitem{Ref11}
B. Jourdain, S. Méléard, W. A. Woyczynski, L\'evy flights in evolutionary ecology. J. Math. Biol. \textbf{65}, 677–707 (2012)

\bibitem{Ref12}
E. R. Weeks, T. H. Solomon, J. S. Urbach, H. L. Swinney, Observation of anomalous diffusion andL\'evy flights BT -L\'evy Flights and Related Topics in Physics in M. F. Shlesinger, G. M. Zaslavsky, U. Frisch, Eds. Springer Berlin Heidelberg (1995)

\bibitem{Ref13}
B. Böttcher, Feller Processes: The Next Generation in Modeling. Brownian Motion,L\'evy Processes and Beyond. PLoS One \textbf{5}, e15102 (2010)

\bibitem{Ref14}
Y. Zheng, F. Yang, J. Duan, X. Sun, L. Fu, J. Kurths, The maximum likelihood climate change for global warming under the influence of greenhouse effect andL\'evy noise. Chaos \textbf{30}, 013132 (2020)

\bibitem{Ref15}
R. Cai, X. Chen, J. Duan, J. Kurths, X. Li, L\'evy noise-induced escape in an excitable system. J. Stat. Mech. Theory Exp. \textbf{6}, 063503 (2017)

\bibitem{Ref16}
R. Cai, Z. He, Y. Liu, J. Duan, J. Kurths, X. Li, Effects of  L\'evy noise on the Fitzhugh-Nagumo model: A perspective on the maximal likely trajectories. J. Theor. Biol. \textbf{480}, 166–174 (2019)

\bibitem{Ref17}
A. Patel, B. Kosko, Stochastic Resonance in Continuous and Spiking Neuron Models with L\'evy Noise. IEEE Trans. Neural Networks \textbf{19}, 1993–2008 (2008)

\bibitem{Ref18}
A. La Cognata, D. Valenti, A. A. Dubkov, B. Spagnolo, Dynamics of two competing species in the presence of  L\'evy noise sources. Phys. Rev. E \textbf{82}, 11121 (2010)

\bibitem{Ref19}
X. Cheng, H. Wang, X. Wang, J. Duan, X. Li, Most probable transition pathways and maximal likely trajectories in a genetic regulatory system. Phys. A Stat. Mech. its Appl. \textbf{531}, 121779 (2019)

\bibitem{Ref20}
L. Serdukova, Y. Zheng, J. Duan, J. Kurths, Stochastic basins of attraction for metastable states. Chaos \textbf{26}, 073117 (2016)

\bibitem{Ref21}
F. Wu, X. Chen, Y. Zheng, J. Duan, J. Kurths, X. Li, L\'evy noise induced transition and enhanced stability in a gene regulatory network. Chaos \textbf{28}, 075510 (2018)

\bibitem{Ref22}
C. Guarcello, D. Valenti, A. Carollo, B. Spagnolo, Effects of  L\'evy noise on the dynamics of sine-Gordon solitons in long Josephson junctions. J. Stat. Mech. Theory Exp. \textbf{2016}, 054012 (2016)

\bibitem{Ref23}
Q. Liu, D. Jiang, T. Hayat, B. Ahmad, Analysis of a delayed vaccinated SIR epidemic model with temporary immunity andL\'evy jumps. Nonlinear Anal. Hybrid Syst. \textbf{27}, 29–43 (2018)

\bibitem{Ref24}
Y. Li, J. Duan, X. Liu, Y. Zhang, Most Probable Dynamics of Stochastic Dynamical Systems with Exponentially Light Jump Fluctuations. Chaos \textbf{30}, 063142 (2020)

\bibitem{Ref25}
Y. Li, J. Duan, A Data-Driven Approach for Discovering Stochastic Dynamical Systems with Non-Gaussian Levy Noise. Phys. D Nonlinear Phenom. \textbf{417}, 132830 (2021)

\bibitem{Ref26}
R. Cont, P. Tankov, \emph{Financial Modelling with Jump Processes}, CRC Press, UK (2003)

\bibitem{Ref27}
A. Janicki, A. Weron, \emph{Simulation and Chaotic Behavior of $\alpha$-stable Stochastic Processes}, Wroclaw University of Technology, Hugo Steinhaus Center (1994)

\bibitem{Ref28}
X. Sun, X. Li, Y. Zheng, Governing equations for probability densities of Marcus stochastic differential equations with L\'evy noise. Stochastics Dyn. \textbf{17}, 1750033 (2016)

\bibitem{Ref29}
J. Duan, \emph{An Introduction to Stochastic Dynamics}, Cambridge University Press, New York (2015)

\bibitem{Ref30}
Q. Liu, Y. Jia, Fluctuations-induced switch in the gene transcriptional regulatory system. Phys. Rev. E \textbf{70}, 041907 (2004)

\bibitem{Ref31}
R. Adler, R. Feldman, M. Taqqu, \emph{A Practical Guide to Heavy Tails: Statistical Techniques and Applications}, Birkhäuser Boston (1998)

\bibitem{Ref32}
A. E. Kyprianou, \emph{Fluctuations of  L\'evy Processes with Applications: Introductory Lectures}, Springer Berlin Heidelberg (2014)

\bibitem{Ref33}
O. E. Barndorff-Nielsen, T. Mikosch, S. I. Resnick, \emph{L\'evy Processes: Theory and Applications}, Springer Science and Business Media, Birkhäuser Boston (2001)

\bibitem{Ref34}
D. Applebaum, \emph{L\'evy Processes and Stochastic Calculus}, 2nd ed, Cambridge University Press, New York (2009)

\bibitem{Ref35}
S. G. Krantz, H. R. Parks, \emph{Geometric Integration Theory}, Springer Science and Business Media, Birkhäuser Boston (2008)

\bibitem{Ref36}
Y. Li, \url{https://github.com/liyangnuaa/Machine-learning-for-asymmetric-Levy-motion}. GitHub (2020)

\end{thebibliography}


\end{document}